# On a third order CWENO boundary treatment with application to networks of hyperbolic conservation laws

Alexander Naumann*, Oliver Kolb*, Matteo Semplice†

**Abstract**

High order numerical methods for networks of hyperbolic conservation laws have recently gained increasing popularity. Here, the crucial part is to treat the boundaries of the single (one-dimensional) computational domains in such a way that the desired convergence rate is achieved in the smooth case but also stability criterions are fulfilled, in particular in the presence of discontinuities. Most of the recently proposed methods rely on a WENO extrapolation technique introduced by Tan and Shu in [*J. Comput. Phys.* 229, pp. 8144–8166 (2010)]. Within this work, we refine and in a sense generalize these results for the case of a third order scheme. Numerical evidence for the analytically found parameter bounds is given as well as results for a complete third order scheme based on the proposed boundary treatment.

## 1 Introduction

There are several applications where networks of hyperbolic conservation laws are involved, e.g., gas and water supply networks [6, 12, 23], traffic flow [5, 9, 15, 17] or blood flow [26, 27]. Those models mainly consist of one-dimensional domains on the edges of an underlying graph complemented with appropriate coupling and boundary conditions at the nodes of the network. Here, the numerical treatment of the boundaries of the single computational domains requires particular attention, since in a network, information leaving one edge may enter another so that even solution components which may be neglected in other applications have to be treated with great care.

The numerous recent publications on higher order boundary treatment show the growing interest in this topic, as well as the potential of high order methods for networks of hyperbolic conservation laws. A second order boundary treatment for networks was recently proposed in [2]. Higher order methods based on WENO and ADER schemes have been considered in [3, 4, 10, 26]. The latter all rely on WENO reconstruction techniques within the boundary cells of each computational domain, mainly based on [30].

*University of Mannheim, Department of Mathematics, 68131 Mannheim, Germany (anaumann@mail.uni-mannheim.de, kolb@uni-mannheim.de).
†Università di Torino, Dipartmento di Matematica, 20123 Torino, Italy (matteo.semplice@unito.it).

There are several recent publications on the parameter choice within WENO schemes (cf. [1,11,13, 20–22,28]) or adapted approaches (cf. [8,14,18]) to avoid order reduction or delayed convergence in the classical WENO schemes. In this work, we provide the development and rigorous analysis of a third order reconstruction procedure for boundary cells, which is based on the CWENO3 reconstruction procedure introduced in [24]. We refine and in a sense generalize the results of [30] for the case of a third order scheme. The considered spatial reconstruction procedure can for instance be applied within the mentioned ADER approaches [3,4,10,26]. Nevertheless, within this work, we will use a third order TVD Runge-Kutta method [16] for the time integration and the presented boundary treatment is coupled with the well-known third order central CWENO3 scheme of [24], with parameter settings according to [20] in the interior of the computational domain.

The outline of this paper is as follows. Section 2 contains the proposed CWENO reconstruction procedure for boundary cells based on cell averages (Section 2.1) and the main analytical results (Section 2.2). Sufficient conditions that ensure the applicability of the convergence results are considered in Section 3. Section 4 is devoted to numerical results. In particular we give some numerical evidence of the performed analysis in Section 4.1, and consider several test scenarios for a full scheme in Sections 4.2 (shock-acoustic interaction), 4.3 (small traffic network), 4.4 (time-dependent boundary conditions), 4.5 (open channel flow) and 4.6 (dam break against a wall). Conclusions as well as an outlook to future work are made in Section 5.

# 2 Numerical scheme and convergence results

In the following section, we develop a third order CWENO reconstruction procedure for boundary cells of a one-dimensional computational domain based on cell averages. Since at the boundary only one-sided information is available, we refer to the ideas of the WENO extrapolation proposed in [30]. As in the interior of the domain, we stick by the design of central WENO schemes as in [24] and additionally consider mesh size dependent parameters within our scheme as in [20]. Convergence results for the proposed reconstruction procedure in different situations are given in Section 2.2. Feasible parameter regions to fulfill the assumptions of these convergence results will be analysed in Section 3.

## 2.1 Reconstruction procedure

The presented reconstruction procedure is independent of the time variable and it therefore suffices to consider functions $u = u(x)$. Further, we assume an equidistant grid with mesh size $h$, grid points $x_j = x_0 + jh$ and finite volumes

$$I_j = [x_j - \frac{h}{2}, x_j + \frac{h}{2}] = [x_{j-0.5}, x_{j+0.5}], \qquad j \in \{0.5, \dots, N-0.5\},$$

as illustrated in Figure 1. Note that the cell interfaces correspond to integer values of $j$ and the cell centers to half-integer values.

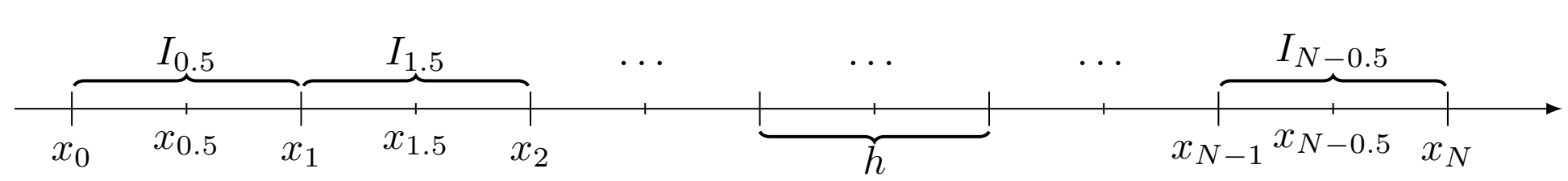

Figure 1: Illustration of the grid notation.

Based on the cell averages over $I_j$,

$$\bar{u}_j = \frac{1}{h} \int_{x_{j-0.5}}^{x_{j+0.5}} u(x)dx,$$

the aim is to reconstruct the underlying function $u$ by a piecewise polynomial approximation. In the interior cells we employ the third order CWENO3 reconstruction with $h$-dependent parameters, which has similar accuracy and stability properties as the proposed scheme, see [20]. Here, we are only interested in the reconstruction within the boundary cells $I_{0.5}$ and $I_{N-0.5}$ and will only consider the boundary treatment at the left boundary of the computational domain. The treatment of the right boundary is symmetrical.

For the reconstruction at the left boundary, we will utilize the cell averages $\bar{u}_{0.5}$, $\bar{u}_{1.5}$ and $\bar{u}_{2.5}$, and combine the ideas of the WENO extrapolation in [30], the weight design of CWENO schemes of [24] and appropriate mesh size dependent parameters as in [20]. More precisely, we consider a convex combination of three polynomials $P_0(x)$, $P_1(x)$ and $P_2(x)$ of the form

$$P(x) = \omega_0 P_0(x) + \omega_1 P_1(x) + \omega_2 P_2(x) \tag{1}$$

with weights $\omega_i \geq 0$ for $i \in \{0, 1, 2\}$ and $\omega_0 + \omega_1 + \omega_2 = 1$.

The first polynomial $P_0(x)$ is the constant polynomial that preserves the cell average $\bar{u}_{0.5}$ over $I_{0.5}$, i.e.,

$$P_0(x) = \bar{u}_{0.5}.$$

The unique linear polynomial that preserves the cell average $\bar{u}_{0.5}$ over $I_{0.5}$ and $\bar{u}_{1.5}$ over $I_{1.5}$ yields the second polynomial:

$$P_1(x) = \bar{u}_{1.5} + \frac{\bar{u}_{1.5} - \bar{u}_{0.5}}{h}(x - x_{1.5}).$$

The third polynomial $P_2(x)$ is constructed by utilizing the unique parabola $P_{opt}(x)$ that satisfies

$$\frac{1}{h} \int_{I_{0.5+l}} P_{opt}(x)dx = \bar{u}_{0.5+l} \qquad \forall l \in \{0, 1, 2\}, \tag{2}$$

i.e., the polynomial that preserves the cell averages $\bar{u}_{0.5}$, $\bar{u}_{1.5}$, $\bar{u}_{2.5}$ over $I_{0.5}$, $I_{1.5}$, $I_{2.5}$, respectively. More explicitly, the unique parabola fulfilling (2) is

$$P_{opt}(x) = u_{1.5} + u'_{1.5}(x - x_{1.5}) + \frac{1}{2}u''_{1.5}(x - x_{1.5})^2$$

with

$$u_{1.5} = \bar{u}_{1.5} - \frac{1}{24}(\bar{u}_{2.5} - 2\bar{u}_{1.5} + \bar{u}_{0.5}), \qquad u'_{1.5} = \frac{\bar{u}_{2.5} - \bar{u}_{0.5}}{2h} \qquad \text{and} \qquad u''_{1.5} = \frac{\bar{u}_{2.5} - 2\bar{u}_{1.5} + \bar{u}_{0.5}}{h^2}.$$

Then, we construct the parabola $P_2(x)$ according to

$$P_{opt}(x) = c_0 P_0(x) + c_1 P_1(x) + c_2 P_2(x)$$

with weights $c_i > 0$ for $i \in \{0, 1, 2\}$ and $c_0 + c_1 + c_2 = 1$, and obtain

$$\begin{aligned} P_2(x) = \frac{1}{c_2}\Big( &- c_0\bar{u}_{0.5} + (1 - c_1)\bar{u}_{1.5} - \frac{1}{24}(\bar{u}_{2.5} - 2\bar{u}_{1.5} + \bar{u}_{0.5}) \\ &+ \frac{\bar{u}_{2.5} - 2c_1\bar{u}_{1.5} + (2c_1 - 1)\bar{u}_{0.5}}{2h}(x - x_{1.5}) + \frac{\bar{u}_{2.5} - 2\bar{u}_{1.5} + \bar{u}_{0.5}}{2h^2}(x - x_{1.5})^2\Big). \end{aligned} \tag{3}$$

Rewriting the polynomials $P_i(x)$ by substituting

$$\sigma_r = \bar{u}_{2.5} - \bar{u}_{1.5} \qquad \text{and} \qquad \sigma_l = \bar{u}_{1.5} - \bar{u}_{0.5}$$

yields

$$\begin{aligned}
P_0(x) &= \bar{u}_{0.5}, \\
P_1(x) &= \bar{u}_{1.5} + \frac{\sigma_l}{h}(x - x_{1.5}), \\
P_2(x) &= \frac{1}{c_2}\Big( - c_0\bar{u}_{0.5} + (1-c_1)\bar{u}_{1.5} - \frac{1}{24}(\sigma_r - \sigma_l) \\
&\qquad + \frac{\sigma_r + \sigma_l - 2c_1\sigma_l}{2h}(x - x_{1.5}) + \frac{\sigma_r - \sigma_l}{2h^2}(x - x_{1.5})^2\Big).
\end{aligned} \tag{4}$$

The weights $\omega_i$ in (1) are chosen in a WENO fashion, according to [25],

$$\omega_i = \frac{\alpha_i}{\sum\limits_{k=0}^{2} \alpha_k} \qquad \text{where} \qquad \alpha_i = \frac{c_i}{(\epsilon(h) + IS_i)^p} \tag{5}$$

with smoothness indicators $IS_i$. The importance of $\epsilon$ within the computation of the $\alpha_i$ has changed from its initial purpose to merely prevent division by zero to a crucial parameter that affects the accuracy and stability of the entire scheme. The originally constant choice of $\epsilon$ was replaced by an $h$-dependent choice $\epsilon(h)$ in several works, for instance [1, 11, 20, 28], since superior convergence properties were observed. Therefore, here we also utilize an $h$-dependent

$$\epsilon(h) = Kh^q \tag{6}$$

with $K \in \mathbb{R}_{>0}$ and $q \in \mathbb{R}_{\geq 0}$.

The smoothness indicators $IS_i$ are computed as proposed in [19], i.e.,

$$IS_i = \sum_{k=1}^{2} \int_{x_0}^{x_1} h^{2k-1} \left(P_i^{(k)}(x)\right)^2 dx. \tag{7}$$

Note that the smoothness in the outermost interval $I_{0.5}$ is investigated. Explicitly computing the smoothness indicators for the polynomials $P_i(x)$ of equations (4) yields

$$IS_0 = 0, \qquad IS_1 = \sigma_l^2, \qquad IS_2 = \frac{1}{c_2^2}\left(\tfrac{4}{3}\sigma_r^2 + (c_1 - \tfrac{11}{3})\sigma_r\sigma_l + (\tfrac{10}{3} - 3c_1 + c_1^2)\sigma_l^2\right). \tag{8}$$

Since the polynomials $P_i(x)$ are affected by the value of $u$ from different regions, we introduce the stencils

$$S_0 = I_{0.5}, \qquad S_1 = I_{0.5} \cup I_{1.5} \qquad \text{and} \qquad S_2 = I_{0.5} \cup I_{1.5} \cup I_{2.5}$$

for later use.

Motivated by the $h$-dependent choice $\epsilon(h)$, we will also consider an $h$-dependent choice for the optimal weights $c_i$,

$$c = [c_0, c_1, c_2] = [K_0 h^{\gamma_0},\, K_1 h^{\gamma_1},\, 1 - K_0 h^{\gamma_0} - K_1 h^{\gamma_1}] \tag{9}$$

with $K_0, K_1 \in \mathbb{R}_{>0}$ and $\gamma_0, \gamma_1 \in \mathbb{R}_{\geq 0}$. Further, we require that

$$\gamma_0 \geq \gamma_1 \qquad \text{and} \qquad \gamma_0 \geq 1,$$

since this will turn out to be necessary in one of the considered cases anyway and additionally simplifies the analysis of other cases. Note that the presence of $P_0$ in (1) could degrade the accuracy to first order unless its coefficient $\omega_0$ tended to zero fast enough in case of smoothness. Intuitively, the second of the above requirements is designed to tackle this and the first one ensures that, if the data are smooth enough, the information brought by $P_1$ (which are second order accurate) have a greater influence in the final reconstructed polynomial.

**Remark 2.1.** *Note that one has to make sure that the resulting weight $c_2$ is positive. Therefore, a robust implementation should for instance apply*

$$c_i = \min\{K_i h^{\gamma_i}, 0.25\}$$

*for $i \in \{0, 1\}$, which does not affect the results for the limit case $h \to 0$ below.*

**Remark 2.2.** *We remark that the well-known linear scaling limiter from [31, 32] can be directly integrated in the presented boundary treatment to achieve a maximum-principle-satisfying and positivity-preserving third order scheme if desired.*

## 2.2 Main results

We now present the main results regarding the convergence of our reconstruction procedure. Following the considerations in [30], we examine three different cases: the smooth case, a single discontinuity in the interval $I_{2.5}$ and a single discontinuity in the interval $I_{1.5}$.

### 2.2.1 Smooth Case

We start by considering the smooth case, i.e., $u$ is smooth within $S_2$. Here, we will show that the presented reconstruction procedure is third order accurate at every point of $I_{0.5}$.

**Theorem 2.1.** *Let $P(x)$ be given by* (1) *with polynomials* (4) *and optimal weights $c_i$ according to* (9) *with $\gamma_0 \geq 1$. Then, if $u$ is smooth enough in the stencil $S_2$ and the weights $\omega_i$ satisfy*

$$\omega_0 = c_0 + O(h^2), \qquad \omega_1 = c_1 + O(h) \qquad \textit{and} \qquad \omega_2 = c_2 + O(h),$$

*we have*

$$u(x) - P(x) = O(h^3) \qquad \forall x \in S_2.$$

*Proof.* Exploiting that $P_{opt}(x)$ was built in such a way that it is a third order approximation to $u(x)$ in the smooth case, we get

$$u(x) - P(x) = P_{opt}(x) + O(h^3) - P(x).$$

Thus, to conclude the proof, it is sufficient to show that

$$P_{opt}(x) - P(x) = O(h^3).$$

By definition we have

$$P_{opt}(x) - P(x) = \sum_i c_i P_i(x) - \sum_i \omega_i P_i(x) = \sum_i (c_i - \omega_i) P_i(x).$$

Making use of $\sum_i c_i = 1 = \sum_i \omega_i$, we obtain

$$\begin{aligned} P_{opt}(x) - P(x) &= \sum_i (c_i - \omega_i)(P_i(x) - u(x)) \\ &= \underbrace{(c_0 - \omega_0)}_{=O(h^2)} \underbrace{(P_0(x) - u(x))}_{=O(h)} + \underbrace{(c_1 - \omega_1)}_{=O(h)} \underbrace{(P_1(x) - u(x))}_{=O(h^2)} + \underbrace{(c_2 - \omega_2)}_{=O(h)} \underbrace{(P_2(x) - u(x))}_{=O(h^2)} \\ &= O(h^3). \end{aligned}$$

Here, $P_2(x)$ is second order accurate since $\gamma_0 \geq 1$. □

**Remark 2.3.** *One key difference between the set-up presented here and the WENO extrapolation employed in [30] is the definition of $P_2(x)$. In [30] $P_2(x) = P_{opt}(x)$ is used and therewith one gets less flexibility in the choice of the weights $c_i$. Considering*

$$u(x) - P(x) = \sum_i \omega_i u(x) - \sum_i \omega_i P_i(x) = \sum_i \omega_i (u(x) - P_i(x)),$$

*we see that $\omega_i = O(h^{2-i})$ is needed for the approximation to be third order accurate. Instead of focussing on the approximation of the optimal weights $c_i$ by $\omega_i$, the set-up in [30] always requires two $h$-dependent weights $c_0$ and $c_1$ with $\gamma_0 \geq 2$ and $\gamma_1 \geq 1$.*

### 2.2.2 Discontinuity in $I_{2.5}$

Assuming a discontinuity in $I_{2.5}$ and smoothness of $u$ in $S_1$, we next show that our reconstruction procedure is, in this case, second order accurate at every point of $I_{0.5}$.

**Theorem 2.2.** *Let $P(x)$ be given by* (1) *with polynomials* (4)*. Moreover, let $u$ be smooth in the stencil $S_1$ and contain a discontinuity in $I_{2.5}$. If the weights $\omega_i$ satisfy the conditions*

$$\omega_0 = O(h), \qquad \omega_1 = O(1) \qquad \textit{and} \qquad \omega_2 = O(h^2),$$

*then*

$$u(x) - P(x) = O(h^2) \qquad \forall x \in S_1,$$

*and even up to the position of the discontinuity (if $u$ is sufficiently smooth).*

*Proof.* Since $\sum_i \omega_i = 1$, we consider

$$u(x) - P(x) = \sum_i \omega_i u(x) - \sum_i \omega_i P_i(x) = \sum_i \omega_i (u(x) - P_i(x)).$$

Hence, it is sufficient to show $\omega_i(u(x) - P_i(x)) = O(h^2)$ for each $i \in \{0, 1, 2\}$. Making use of the theory of polynomial interpolation and our assumptions on the weights $\omega_i$, we directly get

$$\underbrace{\omega_0}_{=O(h)} \underbrace{(u(x) - P_0(x))}_{=O(h)} = O(h^2),$$
$$\underbrace{\omega_1}_{=O(1)} \underbrace{(u(x) - P_1(x))}_{=O(h^2)} = O(h^2),$$
$$\underbrace{\omega_2}_{=O(h^2)} \underbrace{(u(x) - P_2(x))}_{=O(1)} = O(h^2),$$

and thus $u(x) - P(x) = O(h^2)$. □

### 2.2.3 Discontinuity in $I_{1.5}$

Finally, we consider the case where $u$ contains a discontinuity in $I_{1.5}$ and is smooth in $S_0$. Here, we show that our reconstruction procedure is first order accurate at every point of $I_{0.5}$. Note that in this situation the only reliable data is the cell average $\bar{u}_{0.5}$ so that we cannot aim at better accuracy, and indeed in this difficult situation, the task of the reconstruction procedure is not to cause spurious oscillations that would prevent even first order convergence rate.

**Theorem 2.3.** *Let $P(x)$ be given by* (1) *with polynomials* (4)*. Moreover, let $u$ be smooth in $S_0$ and contain a discontinuity in $I_{1.5}$. If the weights $\omega_i$ satisfy the conditions*

$$\omega_0 = O(1), \qquad \omega_1 = O(h) \qquad \text{and} \qquad \omega_2 = O(h),$$

*then*

$$u(x) - P(x) = O(h) \qquad \forall x \in S_0,$$

*and even up to the position of the discontinuity (if $u$ is sufficiently smooth).*

*Proof.* As above, we consider

$$u(x) - P(x) = \sum_i \omega_i u(x) - \sum_i \omega_i P_i(x) = \sum_i \omega_i (u(x) - P_i(x)).$$

Thus, it is sufficient to show $\omega_i(u(x) - P_i(x)) = O(h)$ for each $i \in \{0, 1, 2\}$. Next, the theory of polynomial interpolation and our assumptions on the weights $\omega_i$ yield

$$\begin{aligned}
\underbrace{\omega_0}_{=O(1)} \underbrace{(u(x) - P_0(x))}_{=O(h)} &= O(h), \\
\underbrace{\omega_1}_{=O(h)} \underbrace{(u(x) - P_1(x))}_{=O(1)} &= O(h), \\
\underbrace{\omega_2}_{=O(h)} \underbrace{(u(x) - P_2(x))}_{=O(1)} &= O(h),
\end{aligned}$$

and consequently $u(x) - P(x) = O(h)$ as desired. □

# 3 Accuracy analysis

The accuracy analysis consists of two parts. First, we will investigate properties of the smoothness indicators in Section 3.1. Subsequently, in Section 3.2, specific conditions for the parameter choice within our scheme that guarantee the applicability of our convergence results from Section 2.2 will be deduced.

As already mentioned above, we consider the restrictions

$$\gamma_0 \geq \gamma_1 \qquad \text{and} \qquad \gamma_0 \geq 1,$$

since this will turn out to be necessary in one of the considered cases and additionally simplify the analysis of other cases. Actually, we will even claim $\gamma_0 \geq 1 + \gamma_1$ below to achieve all desired properties, but only $\gamma_0 \geq \gamma_1$ is necessary for the presented intermediate results to be valid. Furthermore, we impose a restriction on $q$ based on the results of [20, 21], namely

$$q \leq 2,$$

since one can expect similar difficulties for the boundary treatment in the case $q > 2$.

## 3.1 The smoothness indicators

The smoothness indicators $IS_i$ play a crucial role in the computation of the weights $\omega_i$ of our reconstruction. We will now take a closer look at their value in the smooth case.

**Lemma 3.1.** *Assuming $u$ is sufficiently smooth in the respective stencil, the smoothness indicators $IS_i$ can be written in the form*

$$\begin{aligned} IS_0 &= 0, \\ IS_1 &= \bar{a}h^2 - \bar{b}h^3 + O(h^4), \\ IS_2 &= \bar{c}_1\bar{a}h^2 + \bar{c}_2\bar{b}h^3 + O(h^4), \end{aligned}$$

*where $\bar{a}, \bar{b} \in \mathbb{R}$ do not depend on the mesh size $h$, and $\bar{c}_1, \bar{c}_2 \in \mathbb{R}$ only depend on the choice of the optimal weights $c_i$.*

*Proof.* Obviously, $IS_0 = 0$ holds and we focus on the analysis of $IS_1$ and $IS_2$. Since the smoothness indicators are based on cell averages, it is convenient to consider the function

$$\bar{u}(x) = \frac{1}{h}\int_{x-\frac{h}{2}}^{x+\frac{h}{2}} u(y)\,dy = \frac{1}{h}\Big(U(x+\frac{h}{2}) - U(x-\frac{h}{2})\Big),$$

where $U$ is a primitive function of $u$. Now, for sufficiently smooth $u$, Taylor expansion gives

$$u^{(k-1)}(x_j+d) = u^{(k-1)}(x_j) + d\,u^{(k)}(x_j) + \frac{d^2}{2}u^{(k+1)}(x_j) + \frac{d^3}{3!}u^{(k+2)}(x_j) + \dots$$

and therewith

$$\begin{aligned} \bar{u}^{(k)}(x_j) = \frac{1}{h}\Big(u^{(k-1)}(x_j+\frac{h}{2}) - u^{(k-1)}(x_j-\frac{h}{2})\Big) \\ = \frac{h^0}{1!\cdot 2^0}u^{(k)}(x_j) + \frac{h^2}{3!\cdot 2^2}u^{(k+2)}(x_j) + \frac{h^4}{5!\cdot 2^4}u^{(k+4)}(x_j) + O(h^6). \end{aligned}$$

Next, the Taylor expansion of $\bar{u}(x_j+h) = \bar{u}_{j+1}$ can be written as

$$\bar{u}_{j+1} = \bar{u}_j + \sum_{k=1}^{3}\frac{h^k}{k!}\bar{u}^{(k)}(x_j) + O(h^4) = \bar{u}_j + a_1h + a_2h^2 + a_3h^3 + O(h^4),$$

where

$$a_n = u^{(n)}(x_j)\sum_{k=0}^{\lceil\frac{n}{2}\rceil-1}\frac{1}{(n-2k)!}\frac{1}{(2k+1)!\cdot 2^{2k}}.$$

Analogously, the Taylor expansion of $\bar{u}(x_j-h) = \bar{u}_{j-1}$ yields

$$\bar{u}_{j-1} = \bar{u}_j + (-1)a_1h + (-1)^2a_2h^2 + (-1)^3a_3h^3 + O(h^4).$$

We define

$$\bar{a} = a_1^2, \qquad \bar{b} = 2a_1a_2, \qquad \bar{c}_1 = \frac{(1-c_1)^2}{c_2^2} \qquad \text{and} \qquad \bar{c}_2 = \frac{-c_1^2+3c_1-2}{c_2^2}.$$

Then, the smoothness indicators $IS_1$ and $IS_2$ can be written in the desired form based on the expansions of $\bar{u}_{j-1}$ and $\bar{u}_{j+1}$ with $j = 1.5$,

$$IS_1 = \sigma_l^2 = (\bar{u}_{1.5} - \bar{u}_{0.5})^2 = a_1^2h^2 - 2a_1a_2h^3 + O(h^4) = \bar{a}h^2 - \bar{b}h^3 + O(h^4),$$

and

$$IS_2 = \bar{c}_1\bar{a}h^2 + \bar{c}_2\bar{b}h^3 + O(h^4).$$

□

**Remark 3.1.** *For later use, we remark that $\bar{c}_1 = 1 + O(h^{\gamma_0})$.*

Based on the given expansions of the smoothness indicators, we can deduce further properties which will become of particular importance in the analysis of the weights $\omega_i$.

**Lemma 3.2.** *Assuming $u$ is sufficiently smooth in $S_2$, we can write*

$$\frac{IS_i - IS_k}{\epsilon(h) + IS_k} = e_{ik}h^{2-q} + O(h^{3-q}) + e_{ik}O(h^{2-q+\delta}),$$

*where all $e_{ik}$ only depend on the choice of the optimal weights $c_i$ and*

$$\delta = \begin{cases} 1 & \text{if } q = 2, \\ 2-q & \text{if } q < 2. \end{cases}$$

*Furthermore,*

$$\frac{\epsilon(h) + IS_i}{\epsilon(h) + IS_k} = 1 + O(h^{2-q})$$

*holds.*

*Proof.* Using the results from Lemma 3.1 for the smoothness indicators $IS_i$, we can calculate the differences

$$\begin{aligned} IS_1 - IS_2 &= (1-\bar{c}_1)\bar{a}h^2 - (1+\bar{c}_2)\bar{b}h^3 + O(h^4), \\ IS_1 - IS_0 &= \bar{a}h^2 - \bar{b}h^3 + O(h^4), \\ IS_2 - IS_0 &= \bar{c}_1\bar{a}h^2 + \bar{c}_2\bar{b}h^3 + O(h^4). \end{aligned} \tag{10}$$

Since $IS_i - IS_k = (-1)(IS_k - IS_i)$, we have

$$IS_i - IS_k = c_{ik}h^2 + O(h^3)$$

in all cases for appropriately chosen $c_{ik}$, which are functions of the optimal weights $c_j$. Moreover, by exploiting our assumption that $q \leq 2$, we have

$$\begin{aligned} \epsilon(h) + IS_0 &= Kh^q, \\ \epsilon(h) + IS_1 &= h^q(K + \bar{a}h^{2-q}) + O(h^3), \\ \epsilon(h) + IS_2 &= h^q(K + \bar{c}_1\bar{a}h^{2-q}) + O(h^3), \end{aligned}$$

which can be summarized to $\epsilon(h) + IS_k = \tilde{K}_k h^q + O(h^{q+\delta})$ with

$$\delta = \begin{cases} 1 & \text{if } q = 2, \\ 2-q & \text{if } q < 2. \end{cases}$$

Altogether,

$$\frac{IS_i - IS_k}{\epsilon(h) + IS_k} = \frac{c_{ik}h^2 + O(h^3)}{\tilde{K}_k h^q + O(h^{q+\delta})} = e_{ik}h^{2-q} + O(h^{3-q}) + e_{ik}O(h^{2-q+\delta})$$

with suitable $e_{ik}$, which are functions of the optimal weights $c_j$.

The second result

$$\frac{\epsilon(h) + IS_i}{\epsilon(h) + IS_k} = 1 + O(h^{2-q})$$

directly follows from

$$\frac{\epsilon(h) + IS_k - IS_k + IS_i}{\epsilon(h) + IS_k} = \frac{\epsilon(h) + IS_k}{\epsilon(h) + IS_k} + \frac{IS_i - IS_k}{\epsilon(h) + IS_k} = 1 + \frac{IS_i - IS_k}{\epsilon(h) + IS_k}.$$

□

**Remark 3.2.** *Considering* (10) *and Remark 3.1, we get*

$$e_{12} = O(h^{\gamma_0}) \qquad \text{and accordingly} \qquad e_{21} = O(h^{\gamma_0}).$$

### 3.2 The weights

The weights $\omega_i$ constitute an integral part of our reconstruction. Utilizing the available information on the smoothness indicators, we now derive sufficient conditions for the parameters within our reconstruction procedure that will guarantee the applicability of the convergence results of Section 2.2. Similar to above, we will consider three cases.

#### 3.2.1 Smooth Case

We start by examining the smooth case, where $u$ is smooth in $S_2$.

**Theorem 3.1.** *We assume that u is smooth enough in the stencil $S_2$. Let the weights $c_i$ be given as $[K_0h^{\gamma_0}, K_1h^{\gamma_1}, 1-K_0h^{\gamma_0}-K_1h^{\gamma_1}]$ with $K_0, K_1 \in \mathbb{R}_{>0}$ and $\gamma_0, \gamma_1 \in \mathbb{R}_{\geq 0}$, where $\gamma_0 \geq \gamma_1$ and $\gamma_0 \geq 1$. Further, we assume $p \in \mathbb{N}$ and $q \leq 2$ in the definition of the weights $\omega_i$ in (5). Then*

$$\omega_0 = c_0 + O(h^{\gamma_0+2-q}), \qquad \omega_1 = c_1 + O(h^{\gamma_1+3-q}) \qquad \text{and} \qquad \omega_2 = c_2 + O(h^{l+2-q})$$

*holds with $l = \min\{\gamma_0, \gamma_1 + 1\}$. In particular, if*

$$\gamma_0 \geq q$$

*is fulfilled, then the weights $\omega_i$ fulfill the hypotheses of Theorem 2.1, i.e.,*

$$\omega_0 = c_0 + O(h^2), \qquad \omega_1 = c_1 + O(h) \qquad \text{and} \qquad \omega_2 = c_2 + O(h).$$

*Proof.* Using the generalization of the third binomial formula in combination with the results from Lemma 3.2 (cf. [1]), we obtain

$$\begin{aligned}\frac{1}{(\epsilon(h)+IS_k)^p} &= \frac{1}{(\epsilon(h)+IS_i)^p}\left(1+\underbrace{\frac{IS_i-IS_k}{\epsilon(h)+IS_k}}_{=e_{ik}h^{2-q}+O(h^{3-q})+e_{ik}O(h^{2-q+\delta})}\sum_{t=0}^{p-1}\Bigg(\underbrace{\frac{\epsilon(h)+IS_i}{\epsilon(h)+IS_k}}_{=1+O(h^{2-q})}\Bigg)^t\right)\\ &= \frac{1}{(\epsilon(h)+IS_i)^p}\left(1+pe_{ik}h^{2-q}+O(h^{3-q})+e_{ik}O(h^{2-q+\delta})+e_{ik}O(h^{4-2q})\right).\end{aligned}$$

Defining $z = 2-q+\min(\delta, 2-q)$ and considering the definition of $\alpha_k$ (equation (5)) yields

$$\alpha_k = \frac{c_k}{(\epsilon(h)+IS_k)^p} = \frac{c_k}{(\epsilon(h)+IS_i)^p}\left(1+pe_{ik}h^{2-q}+O(h^{3-q})+e_{ik}O(h^z)\right).$$

Therewith we get for $i \in \{0,1,2\}$

$$\begin{aligned}\sum_k \alpha_k &= \frac{1}{(\epsilon(h)+IS_i)^p}\Big(c_i+\sum_{k\neq i}c_k\left(1+pe_{ik}h^{2-q}+O(h^{3-q})+e_{ik}O(h^z)\right)\Big)\\ &= \frac{1}{(\epsilon(h)+IS_i)^p}\Bigg(\underbrace{\sum_k c_k}_{=1}+\Big(\underbrace{p\sum_{k\neq i}c_ke_{ik}}_{=:-f_i}\Big)h^{2-q}+\sum_{k\neq i}c_kO(h^{3-q})+\sum_{k\neq i}c_ke_{ik}O(h^z)\Bigg)\\ &= \frac{1}{(\epsilon(h)+IS_i)^p}\left(1-pf_ih^{2-q}+O(h^{3-q+r})+O(h^{z+s})\right),\end{aligned}$$

where

$$r = \begin{cases}0 & \text{for } i \in \{0,1\},\\ \gamma_1 & \text{for } i = 2,\end{cases}$$

and

$$s = \begin{cases} 0 & \text{for } i = 0, \\ \gamma_0 & \text{for } i \in \{1, 2\} \end{cases}$$

due to Remark 3.2 and $c_0 = O(h^{\gamma_0})$. For the same reason, we have $f_1 = O(h^{\gamma_0})$ and $f_2 = O(h^{\gamma_0})$.

We can now compute $\omega_i$ as

$$\omega_i = \frac{\alpha_i}{\sum_k \alpha_k} = \frac{c_i}{1 - pf_i h^{2-q} + O(h^{3-q+r}) + O(h^{z+s})}$$
$$= c_i \left(1 + pf_i h^{2-q} + O(h^{3-q+r}) + O(h^{z+s})\right).$$

With the knowledge of $r$, $s$ and the properties of $f_i$ from above, we finally get (together with the assumptions on $\gamma_0$, $\gamma_1$ and $q$, and the fact that $z \geq 2 - q$)

$$\begin{aligned}
\omega_0 &= c_0 + O(h^{\gamma_0+2-q}) + O(h^{\gamma_0+3-q}) + O(h^{\gamma_0+z}) = c_0 + O(h^{\gamma_0+2-q}), \\
\omega_1 &= c_1 + O(h^{\gamma_1+\gamma_0+2-q}) + O(h^{\gamma_1+3-q}) + O(h^{\gamma_1+z+\gamma_0}) = c_1 + O(h^{\gamma_1+3-q}), \\
\omega_2 &= c_2 + O(h^{\gamma_0+2-q}) + O(h^{3-q+\gamma_1}) + O(h^{z+\gamma_0}) = c_2 + O(h^{l+2-q}),
\end{aligned}$$

where $l = \min(\gamma_0, \gamma_1 + 1)$. □

### 3.2.2 Discontinuity in $I_{2.5}$

The second case under investigation assumes $u$ to be smooth in $S_1$ and a single discontinuity in the interval $I_{2.5}$.

**Theorem 3.2.** *Let $u(x)$ be smooth in $S_1$ and contain a discontinuity in $I_{2.5}$ with jump $\Delta u$. Further, let the weights $c_i$ be given as $[K_0 h^{\gamma_0}, K_1 h^{\gamma_1}, 1 - K_0 h^{\gamma_0} - K_1 h^{\gamma_1}]$ with $K_0, K_1 \in \mathbb{R}_{>0}$ and $\gamma_0, \gamma_1 \in \mathbb{R}_{\geq 0}$, where $\gamma_0 \geq \gamma_1$. Assuming $pq \geq \gamma_1$ and $q \leq 2$, the weights $\omega_i$ satisfy*

$$\omega_0 = O(h^{\gamma_0-\gamma_1}), \qquad \omega_1 = O(1) \qquad \text{and} \qquad \omega_2 = O(h^{pq-\gamma_1}).$$

*In particular, if the parameters additionally fulfill the conditions*

$$\gamma_0 \geq 1 + \gamma_1 \qquad \text{and} \qquad pq \geq 2 + \gamma_1,$$

*the weights $\omega_i$ satisfy the hypotheses of Theorem 2.2, i.e.,*

$$\omega_0 = O(h), \qquad \omega_1 = O(1) \qquad \text{and} \qquad \omega_2 = O(h^2).$$

*Proof.* The underlying $u$ is smooth in stencil $S_1$. Thus, we get from Lemma 3.1

$$\begin{aligned}
IS_0 &= 0, \\
IS_1 &= \bar{a}h^2 + \bar{b}h^3 + O(h^4).
\end{aligned}$$

Due to the discontinuity in $I_{2.5}$, we further have

$$IS_2 = \Theta(\Delta u^2).$$

According to (5) we obtain

$$\begin{aligned}
\alpha_0 &= \frac{c_0}{(Kh^q + 0)^p} = \Theta(h^{\gamma_0 - pq}), \\
\alpha_1 &= \frac{c_1}{(h^q(K + \bar{a}h^{2-q} + O(h^{3-q})))^p} = \Theta(h^{\gamma_1 - pq}), \\
\alpha_2 &= \frac{c_2}{(Kh^q + \Theta(\Delta u^2))^p} = \Theta(1).
\end{aligned}$$

Since $\gamma_1 \le pq$ and $\gamma_1 \le \gamma_0$, we have

$$\sum_k \alpha_k = \Theta(h^{\gamma_0 - pq}) + \Theta(h^{\gamma_1 - pq}) + \Theta(1) = \Theta(h^{\gamma_1 - pq}).$$

Finally, we obtain the weights as

$$\omega_0 = O(h^{\gamma_0 - \gamma_1}), \qquad \omega_1 = O(1) \qquad \text{and} \qquad \omega_2 = O(h^{pq - \gamma_1}).$$

□

### 3.2.3 Discontinuity in $I_{1.5}$

In the last case we consider a single discontinuity of $u$ in $I_{1.5}$.

**Theorem 3.3.** *Let $u(x)$ be smooth in $S_0$ and contain a discontinuity in $I_{1.5}$ with jump $\Delta u$. Further, let the weights $c_i$ be given as $[K_0 h^{\gamma_0}, K_1 h^{\gamma_1}, 1 - K_0 h^{\gamma_0} - K_1 h^{\gamma_1}]$ with $K_0, K_1 \in \mathbb{R}_{>0}$ and $\gamma_0, \gamma_1 \in \mathbb{R}_{\ge 0}$. If $pq \ge \gamma_0$, then the weights $\omega_i$ fulfill*

$$\omega_0 = O(1), \qquad \omega_1 = O(h^{\gamma_1 + pq - \gamma_0}) \qquad \textit{and} \qquad \omega_2 = O(h^{pq - \gamma_0}).$$

*In particular, if the parameters fulfill*

$$pq \ge 1 + \gamma_0,$$

*then the weights $\omega_i$ satisfy the hypotheses of Theorem 2.3, i.e.,*

$$\omega_0 = O(1), \qquad \omega_1 = O(h) \qquad \textit{and} \qquad \omega_2 = O(h).$$

*Proof.* Computing the smoothness indicators according to (8) yields

$$IS_0 = 0, \qquad IS_1 = \Theta(\Delta u^2) \qquad \text{and} \qquad IS_2 = \Theta(\Delta u^2).$$

The resulting values for $\alpha$ are

$$\begin{aligned} \alpha_0 &= \frac{c_0}{(Kh^q + 0)^p} = \Theta(h^{\gamma_0 - pq}), \\ \alpha_1 &= \frac{c_1}{(Kh^q + \Theta(\Delta u^2))^p} = \Theta(h^{\gamma_1}), \\ \alpha_2 &= \frac{c_2}{(Kh^q + \Theta(\Delta u^2))^p} = \Theta(1). \end{aligned}$$

Therewith we get

$$\sum_k \alpha_k = \Theta(h^{\gamma_0 - pq}) + \Theta(h^{\gamma_1}) + \Theta(1) = \Theta(h^{\gamma_0 - pq})$$

and finally

$$\omega_0 = O(1), \qquad \omega_1 = O(h^{\gamma_1 + pq - \gamma_0}) \qquad \text{and} \qquad \omega_2 = O(h^{pq - \gamma_0}).$$

□

Summarizing, the following set of conditions is sufficient to guarantee that we meet the requirements of the results from Section 2.2:

$$q \le 2, \qquad \gamma_0 \ge \max(q, 1 + \gamma_1), \qquad pq \ge 1 + \gamma_0, \tag{11}$$

where $pq \ge 2 + \gamma_1$ from Section 3.2.2 follows from the latter two conditions.

# 4 Numerical results

Within this section, we consider several numerical examples. Section 4.1 is aimed at numerically confirming the analytical results of Section 2.2 for the proposed boundary treatment. In Sections 4.2 until 4.6, the numerical solution of conservation laws is considered. Here, for the spatial semi-discretization in the interior of the computational domain, we apply the CWENO3 reconstruction procedure of [24] with the new weights designed in [20]. Finally, a third order TVD Runge-Kutta scheme [16] is used for the time integration.

For completeness we give a brief description of the spatial semi-discretization used in Sections 4.2 until 4.6, which is a finite volume approach. Starting from a hyperbolic conservation law $u_t + f(u)_x = 0$ with given initial conditions $u(x,0) = u_0(x)$, averaging over the grid cells $I_j$ yields

$$\frac{\mathrm{d}}{\mathrm{d}t}\bar{u}_j(t) = -\frac{1}{h}\left(f(u(x_{j+0.5},t)) - f(u(x_{j-0.5},t))\right)$$

with initial conditions

$$\bar{u}_j(0) = \frac{1}{h}\int\limits_{x_{j-0.5}}^{x_{j+0.5}} u_0(x)\,dx.$$

Next, the fluxes $f(u(x_{j\pm0.5},t))$ at the cell boundaries are approximated by a numerical flux function $H_{j\pm0.5}(t) = H(u^-_{j\pm0.5}(t), u^+_{j\pm0.5}(t))$, where $H$ is taken to be the local Lax-Friedrichs flux in this work except Sections 4.3 and 4.4, where an exact Riemann solver is applied. The values $u^\pm_{j\pm0.5}(t)$ are reconstructed values from both sides of each cell interface $x_{j\pm0.5}$ based on the cell averages $\bar{u}_j(t)$ at time $t$. The reconstruction within boundary cells is described in Section 2.1. The applied polynomial reconstruction for cells $I_j$ in the interior of the computational domain is very similar but symmetric: it consists of two linear polynomials

$$P_l(x) = \bar{u}_j + \frac{\bar{u}_j - \bar{u}_{j-1}}{h}(x - x_j), \qquad P_r(x) = \bar{u}_j + \frac{\bar{u}_{j+1} - \bar{u}_j}{h}(x - x_j)$$

and a parabola $P_c$, chosen in such a way that

$$P_{opt}(x) = c_l P_l(x) + c_c P_c(x) + c_r P_r(x),$$

where $P_{opt}$ is the unique parabola that preserves the cell averages $\bar{u}_{j-1}$, $\bar{u}_j$ and $\bar{u}_{j+1}$ in $I_{j-1}$, $I_j$ and $I_{j+1}$, respectively. According to [24] we choose $c_l = c_r = 0.25$ and $c_c = 0.5$. The weight design in the applied polynomial reconstruction

$$P(x) = w_l P_l(x) + w_c P_c(x) + w_r P_r(x)$$

is identical to the one presented in Section 2.1 above. For further details we refer to [20].

Note that within all presented results we apply $p = 2$, since beside its role in the conditions on $pq$, the discussion of the parameter $p$ does not give novel insight compared to the results in [20].

## 4.1 Reconstruction procedure

In this section, we will confirm numerically the analytical convergence results of our reconstruction procedure. The same three cases of Section 2.2 will be investigated.

In the smooth case, we consider the function

$$u_s(x) = \sin(2\pi x)$$

and a mesh with its left boundary at the point $x_0 = 0$. Therefore, the reconstruction is based on the cell averages of $u_s(x)$ in the intervals $I_{0.5} = [0, h]$, $I_{1.5} = [h, 2h]$ and $I_{2.5} = [2h, 3h]$. The resulting polynomial is then evaluated at the left boundary $x = x_0 = 0$ and compared to the exact function value $u_s(0) = 0$.

For the discontinuous cases, the function

$$u_d(x) = \sin(2\pi x) + H(x) \qquad \text{with} \qquad H(x) = \begin{cases} \Delta u & \text{if } x > x^*, \\ 0 & \text{else} \end{cases}$$

with $\Delta u = 0.5$ and $x^* = 0.1$ is considered. To get the right setting for the convergence results in Theorems 2.2 and 2.3, we use $x_0 = x^* - 2.5h$ and $x_0 = x^* - 1.5h$, respectively, and evaluate the resulting polynomials at $x = x^*$.

For all test cases, we will use the following four parameter sets:

$$\begin{aligned}
\Sigma_1 &= \{\epsilon = h,\ c = [h, 0.25, 0.75 - h]\}, \\
\Sigma_2 &= \{\epsilon = h^2,\ c = [h^2, 0.25, 0.75 - h^2]\}, \\
\Sigma_3 &= \{\epsilon = 10^{-3},\ c = [h^2, h, 1 - h^2 - h]\}, \\
\Sigma_4 &= \{\epsilon = 10^{-6},\ c = [h^2, h, 1 - h^2 - h]\},
\end{aligned}$$

and consider $h = 0.25 \cdot 2^{-n}$ throughout all subsections of Section 4.1. The first two sets fulfill all necessary conditions (11) from Section 3, whereas $\Sigma_3$ and $\Sigma_4$ contain the original (and with respect to $c$ feasible) choice for the weights $c_i$ but, also according to the original [30], a constant value for $\epsilon$ (within a reasonable range). In addition to those four parameter sets, we will consider further parameter sets for each test case to numerically validate that a violation of the given sufficient conditions usually results in a reduction of the convergence rate (or no convergence at all).

### 4.1.1 Smooth case

We begin with the smooth case. According to Theorem 3.1

$$\gamma_0 \geq q$$

is a critical condition. We expect to observe third order convergence for all parameter sets $\Sigma_1$ to $\Sigma_4$, and additionally consider

$$\Sigma_{5.1} = \{\epsilon = h^2,\ c = [h, 0.25, 0.75 - h]\},$$

which violates this condition ($q = 2 > \gamma_0 = 1$).

The results for the smooth test case for all parameter sets are summarized in Table 1, where we report on the reconstruction errors $|P(0) - u_s(0)|$ and the convergence rates (c.r.). First of all, as expected from Theorems 2.1 and 3.1, all sets $\Sigma_1$ to $\Sigma_4$ finally yield third order accuracy, whereas only second order is achieved for the extra set of parameters $\Sigma_{5.1}$, which does not fulfill the necessary conditions. Nevertheless the results of the first four sets differ in the number of refinement steps to reach the formal order and with respect to the stability of the observed convergence rates. Concerning the latter point, the mesh size dependent choices for $\epsilon$ in $\Sigma_1$ and $\Sigma_2$ are superior to the other choices. Further, the formal order is observed rather late for $\epsilon = 10^{-6}$. Note that both points mentioned so far have similarly been observed for the CWENO3 reconstruction in the interior of the computational domain in [20]. We wish to point out that in the smooth case the absolute errors on the finest levels are the smallest for the two constant choices of $\epsilon$ in $\Sigma_3$ and $\Sigma_4$. However, we will see later on that these two parameter sets do not yield good results in the discontinuous cases.

Table 1: Reconstruction results for $u_s(x)$ at $x = 0$, $h = 0.25 \cdot 2^{-n}$.

| | $\Sigma_1$ $(\gamma_0, \gamma_1) = (1, 0)$ $\epsilon = h$ | | $\Sigma_2$ $(\gamma_0, \gamma_1) = (2, 0)$ $\epsilon = h^2$ | | $\Sigma_3$ $(\gamma_0, \gamma_1) = (2, 1)$ $\epsilon = 10^{-3}$ | | $\Sigma_4$ $(\gamma_0, \gamma_1) = (2, 1)$ $\epsilon = 10^{-6}$ | | $\Sigma_{5.1}$ $(\gamma_0, \gamma_1) = (1, 0)$ $\epsilon = h^2$ | |
|---|---|---|---|---|---|---|---|---|---|---|
| n | Error | c.r. | Error | c.r. | Error | c.r. | Error | c.r. | Error | c.r. |
| 1 | 3.22e-01 | - | 3.59e-01 | - | 3.73e-01 | - | 3.73e-01 | - | 3.71e-01 | - |
| 2 | 9.03e-02 | 1.8 | 1.70e-01 | 1.1 | 1.92e-01 | 1.0 | 1.94e-01 | 0.9 | 1.92e-01 | 0.9 |
| 3 | 1.06e-02 | 3.1 | 6.03e-02 | 1.5 | 5.92e-02 | 1.7 | 9.79e-02 | 1.0 | 9.62e-02 | 1.0 |
| 4 | 1.04e-03 | 3.3 | 1.39e-02 | 2.1 | 1.09e-03 | 5.8 | 4.90e-02 | 1.0 | 4.73e-02 | 1.0 |
| 5 | 1.10e-04 | 3.3 | 2.21e-03 | 2.7 | 1.36e-05 | 6.3 | 2.45e-02 | 1.0 | 2.28e-02 | 1.1 |
| 6 | 1.24e-05 | 3.1 | 2.96e-04 | 2.9 | 3.40e-06 | 2.0 | 1.04e-02 | 1.2 | 1.06e-02 | 1.1 |
| 7 | 1.47e-06 | 3.1 | 3.76e-05 | 3.0 | 4.54e-07 | 2.9 | 4.94e-04 | 4.4 | 4.68e-03 | 1.2 |
| 8 | 1.78e-07 | 3.0 | 4.73e-06 | 3.0 | 5.75e-08 | 3.0 | 4.30e-06 | 6.8 | 1.89e-03 | 1.3 |
| 9 | 2.20e-08 | 3.0 | 5.91e-07 | 3.0 | 7.21e-09 | 3.0 | 3.21e-08 | 7.1 | 6.82e-04 | 1.5 |
| 10 | 2.73e-09 | 3.0 | 7.40e-08 | 3.0 | 9.02e-10 | 3.0 | 4.34e-10 | 6.2 | 2.19e-04 | 1.6 |
| 11 | 3.40e-10 | 3.0 | 9.24e-09 | 3.0 | 1.13e-10 | 3.0 | 1.04e-10 | 2.1 | 6.39e-05 | 1.8 |
| 12 | 4.24e-11 | 3.0 | 1.16e-09 | 3.0 | 1.41e-11 | 3.0 | 1.39e-11 | 2.9 | 1.74e-05 | 1.9 |
| 13 | 5.29e-12 | 3.0 | 1.44e-10 | 3.0 | 1.76e-12 | 3.0 | 1.76e-12 | 3.0 | 4.56e-06 | 1.9 |
| 14 | 6.61e-13 | 3.0 | 1.81e-11 | 3.0 | 2.20e-13 | 3.0 | 2.20e-13 | 3.0 | 1.17e-06 | 2.0 |

### 4.1.2 Discontinuity in $I_{2.5}$

Next, we consider the discontinuous case with single discontinuity in $I_{2.5}$. Therefore, we will reconstruct $u_d(x)$ based on the cell averages in the intervals $I_{0.5} = [x^* - 2.5h, x^* - 1.5h]$, $I_{1.5} = [x^* - 1.5h, x^* - 0.5h]$ and $I_{2.5} = [x^* - 0.5h, x^* + 0.5h]$, and evaluate $|P(x^*) - u_d(x^*)|$. According to Theorem 3.2

$$\gamma_0 \geq 1 + \gamma_1 \qquad \text{and} \qquad pq \geq 2 + \gamma_1$$

are critical conditions, which are fulfilled by $\Sigma_1$ and $\Sigma_2$. Obviously, any constant choice for $\epsilon$ ($q = 0$) violates the second condition. Nevertheless, an interesting point is whether this comes into effect for practically relevant refinement levels. As additional parameter sets we consider here

$$\Sigma_{5.2} = \{\epsilon = h^2,\ c = [h,\ h,\ 1 - 2h]\}$$

and

$$\Sigma_{6.2} = \{\epsilon = h,\ c = [h^2,\ h,\ 1 - h^2 - h]\},$$

which violate each one of those conditions.

The results for this test case are summarized in Table 2. Both $\Sigma_1$ and $\Sigma_2$ reach the expected formal order of accuracy of two (even though $\Sigma_2$ quite late). As in the smooth case, $\Sigma_1$ reaches the desired order faster but the final absolute errors are similar here. The choices $\Sigma_3$ and $\Sigma_4$ with constant values of $\epsilon$ finally do not converge. From the theoretical point of view this must happen as soon as $\epsilon$ outweighs the smoothness indicators of the polynomials that only depend on regions where the solution is smooth. Nevertheless, the accuracy of the parameter set $\Sigma_4$ is comparable to that of $\Sigma_2$ at certain levels in the considered range of $h$. The sharpness of the necessary conditions in Theorem 3.2 is well demonstrated by the results of the parameter sets $\Sigma_{5.2}$ ($\gamma_0 = 1 \ngeq 1+\gamma_1 = 2$) and $\Sigma_{6.2}$ ($pq = 2 \ngeq 2 + \gamma_1 = 3$), for which an order reduction can be clearly observed.

### 4.1.3 Discontinuity in $I_{1.5}$

In the third case, we consider a single discontinuity in $I_{1.5}$ and therefore the reconstruction of $u_d(x)$ based on the cell averages of the intervals $I_{0.5} = [x^* - 1.5h, x^* - 0.5h]$, $I_{1.5} = [x^* - 0.5h, x^* + 0.5h]$ and $I_{2.5} = [x^* + 0.5h, x^* + 1.5h]$. Again, we are interested in the error $|P(x^*) - u_d(x^*)|$.

Table 2: Reconstruction results for $u_d(x)$ at $x^* = x_{2.5}$, $h = 0.25 \cdot 2^{-n}$.

| | $\Sigma_1$ $(\gamma_0,\gamma_1)=(1,0)$ $\epsilon = h$ | | $\Sigma_2$ $(\gamma_0,\gamma_1)=(2,0)$ $\epsilon = h^2$ | | $\Sigma_3$ $(\gamma_0,\gamma_1)=(2,1)$ $\epsilon = 10^{-3}$ | | $\Sigma_4$ $(\gamma_0,\gamma_1)=(2,1)$ $\epsilon = 10^{-6}$ | | $\Sigma_{5.2}$ $(\gamma_0,\gamma_1)=(1,1)$ $\epsilon = h^2$ | | $\Sigma_{6.2}$ $(\gamma_0,\gamma_1)=(2,1)$ $\epsilon = h$ | |
|---|---|---|---|---|---|---|---|---|---|---|---|---|
| n | Error | c.r. | Error | c.r. | Error | c.r. | Error | c.r. | Error | c.r. | Error | c.r. |
| 1 | 9.99e-01 | - | 1.25e-00 | - | 1.38e-00 | - | 1.38e-00 | - | 1.36e-00 | - | 7.83e-02 | - |
| 2 | 2.05e-01 | 2.3 | 6.04e-01 | 1.0 | 7.29e-01 | 0.9 | 7.43e-01 | 0.9 | 7.33e-01 | 0.9 | 2.07e-01 | -1.4 |
| 3 | 7.73e-03 | 4.7 | 2.63e-01 | 1.2 | 2.44e-01 | 1.6 | 3.55e-01 | 1.1 | 3.51e-01 | 1.1 | 2.26e-01 | -0.1 |
| 4 | 5.66e-03 | 0.4 | 7.98e-02 | 1.7 | 1.37e-02 | 4.2 | 1.69e-01 | 1.1 | 1.69e-01 | 1.1 | 2.00e-01 | 0.2 |
| 5 | 1.49e-03 | 1.9 | 1.32e-02 | 2.6 | 2.96e-02 | -1.1 | 8.22e-02 | 1.0 | 8.21e-02 | 1.0 | 1.51e-01 | 0.4 |
| 6 | 3.77e-04 | 2.0 | 1.55e-03 | 3.1 | 1.76e-02 | 0.7 | 4.04e-02 | 1.0 | 4.04e-02 | 1.0 | 9.81e-02 | 0.6 |
| 7 | 9.52e-05 | 2.0 | 1.46e-04 | 3.4 | 2.02e-02 | -0.2 | 1.91e-02 | 1.1 | 2.00e-02 | 1.0 | 5.70e-02 | 0.8 |
| 8 | 2.39e-05 | 2.0 | 7.04e-06 | 4.4 | 3.22e-02 | -0.7 | 3.93e-03 | 2.3 | 9.96e-03 | 1.0 | 3.09e-02 | 0.9 |
| 9 | 6.00e-06 | 2.0 | 1.82e-06 | 2.0 | 5.46e-02 | -0.8 | 1.14e-04 | 5.1 | 4.97e-03 | 1.0 | 1.61e-02 | 0.9 |
| 10 | 1.50e-06 | 2.0 | 8.93e-07 | 1.0 | 8.81e-02 | -0.7 | 1.68e-06 | 6.1 | 2.48e-03 | 1.0 | 8.24e-03 | 1.0 |
| 11 | 3.76e-07 | 2.0 | 2.78e-07 | 1.7 | 1.28e-01 | -0.5 | 5.84e-07 | 1.5 | 1.24e-03 | 1.0 | 4.17e-03 | 1.0 |
| 12 | 9.40e-08 | 2.0 | 7.61e-08 | 1.9 | 1.67e-01 | -0.4 | 7.18e-07 | -0.3 | 6.20e-04 | 1.0 | 2.09e-03 | 1.0 |
| 13 | 2.35e-08 | 2.0 | 1.99e-08 | 1.9 | 1.97e-01 | -0.2 | 1.20e-06 | -0.7 | 3.10e-04 | 1.0 | 1.05e-03 | 1.0 |
| 14 | 5.88e-09 | 2.0 | 5.07e-09 | 2.0 | 2.16e-01 | -0.1 | 2.29e-06 | -0.9 | 1.55e-04 | 1.0 | 5.26e-04 | 1.0 |

Here, according to Theorem 3.3,

$$pq \geq 1 + \gamma_0$$

is a critical condition. Besides the parameter settings $\Sigma_1$ to $\Sigma_4$, we will consider

$$\Sigma_{5.3} = \{\epsilon = h,\ c = [h^{1.5},\ h^{0.5},\ 1 - h^{1.5} - h^{0.5}]\}.$$

The results for the third test case are summarized in Table 3. Again, the choices $\Sigma_1$ and $\Sigma_2$ yield the expected formal order of accuracy, but now $\Sigma_2$ is slightly better than $\Sigma_1$. Further, as in the previous test case, the constant choices of $\epsilon$ within the parameter combinations $\Sigma_3$ and $\Sigma_4$ do not show convergence. The fifth parameter set $\Sigma_{5.3}$, which does not fulfill the additional condition $pq \geq 1+\gamma_0$ of Theorem 3.3, only shows very slow convergence. Formally, a convergence rate of 0.5 should be reached, when reconsidering the proof of Theorem 2.3 together with the statement about the weights $\omega_i$ in Theorem 3.3 ($\omega_2 = O(h^{0.5})$ for the given parameters instead of $\omega_2 = O(h^1)$).

## 4.2 Shock-acoustic interaction

Next, we consider the shock-acoustic interaction example from [29], i.e., the one-dimensional Euler equations

$$\frac{\partial}{\partial t}\begin{pmatrix}\rho\\ m\\ E\end{pmatrix} + \frac{\partial}{\partial x}\begin{pmatrix}m\\ \rho v^2 + p\\ v(E+p)\end{pmatrix} = 0$$

with pressure law $p = (\gamma-1)(E - \rho v^2/2)$, $\gamma = 1.4$, velocity $v = m/\rho$, and initial conditions

$$\begin{pmatrix}\rho\\ v\\ p\end{pmatrix}(x,0) = \begin{cases}(3.857143,\ 2.629369,\ 10.33333) & \text{for } x < -4,\\ (1+\delta\sin(5x),\ 0,\ 1) & \text{for } x \geq -4,\end{cases}$$

with $\delta = 0.2$. For the numerical simulation we use the interval $[-5,5]$ and compare results with and without an artificial boundary at $x = 0$, i.e., once we consider a single computational domain and once we split the domain into two domains and apply the proposed boundary treatment from both sides at the inner boundary. For the scheme parameters we apply

$$\Sigma_1 = \{\epsilon = h,\ c = [h, 0.25, 0.75 - h]\}, \tag{12}$$
$$\Sigma_2 = \{\epsilon = h^2,\ c = [h^2, 0.25, 0.75 - h^2]\} \tag{13}$$

from above. In the interior of the computational domain, we reconstruct each conservative variable using CWENO3 from [20] with the same parameters $p = 2$ and $\epsilon$ for the weight design. The ordinary differential equation resulting from the spatial semi-discretization is solved with a third order TVD Runge-Kutta method [16].

Table 3: Reconstruction results for $u_d(x)$ at $x^* = x_{1.5}$, $h = 0.25 \cdot 2^{-n}$.

| Σ | $\Sigma_1$ $(\gamma_0,\gamma_1)=(1,0)$ $\epsilon = h$ | | $\Sigma_2$ $(\gamma_0,\gamma_1)=(2,0)$ $\epsilon = h^2$ | | $\Sigma_3$ $(\gamma_0,\gamma_1)=(2,1)$ $\epsilon = 10^{-3}$ | | $\Sigma_4$ $(\gamma_0,\gamma_1)=(2,1)$ $\epsilon = 10^{-6}$ | | $\Sigma_{5.3}$ $(\gamma_0,\gamma_1)=(1.5,0.5)$ $\epsilon = h$ | |
|---|---|---|---|---|---|---|---|---|---|---|
| n | Error | c.r. | Error | c.r. | Error | c.r. | Error | c.r. | Error | c.r. |
| 1 | 7.01e-01 | - | 7.32e-01 | - | 7.40e-01 | - | 7.40e-01 | - | 6.17e-01 | - |
| 2 | 2.47e-01 | 1.5 | 3.43e-01 | 1.1 | 3.55e-01 | 1.1 | 3.56e-01 | 1.1 | 5.08e-02 | 3.6 |
| 3 | 1.08e-02 | 4.5 | 1.58e-01 | 1.1 | 1.58e-01 | 1.2 | 1.70e-01 | 1.1 | 1.62e-01 | -1.7 |
| 4 | 7.67e-02 | -2.8 | 7.60e-02 | 1.1 | 1.90e-03 | 6.4 | 8.23e-02 | 1.0 | 2.12e-01 | -0.4 |
| 5 | 9.65e-02 | -0.3 | 3.80e-02 | 1.0 | 1.60e-01 | -6.4 | 4.05e-02 | 1.0 | 2.24e-01 | -0.1 |
| 6 | 8.66e-02 | 0.2 | 1.93e-02 | 1.0 | 2.29e-01 | -0.5 | 2.00e-02 | 1.0 | 2.27e-01 | -0.0 |
| 7 | 6.44e-02 | 0.4 | 9.76e-03 | 1.0 | 2.45e-01 | -0.1 | 9.96e-03 | 1.0 | 2.24e-01 | 0.0 |
| 8 | 4.16e-02 | 0.6 | 4.92e-03 | 1.0 | 2.49e-01 | -0.0 | 4.91e-03 | 1.0 | 2.19e-01 | 0.0 |
| 9 | 2.42e-02 | 0.8 | 2.47e-03 | 1.0 | 2.50e-01 | -0.0 | 2.22e-03 | 1.1 | 2.10e-01 | 0.1 |
| 10 | 1.31e-02 | 0.9 | 1.24e-03 | 1.0 | 2.50e-01 | -0.0 | 1.88e-04 | 3.6 | 1.99e-01 | 0.1 |
| 11 | 6.87e-03 | 0.9 | 6.20e-04 | 1.0 | 2.50e-01 | -0.0 | 3.57e-03 | -4.2 | 1.84e-01 | 0.1 |
| 12 | 3.51e-03 | 1.0 | 3.10e-04 | 1.0 | 2.50e-01 | -0.0 | 1.57e-02 | -2.1 | 1.66e-01 | 0.1 |
| 13 | 1.78e-03 | 1.0 | 1.55e-04 | 1.0 | 2.50e-01 | -0.0 | 5.37e-02 | -1.8 | 1.46e-01 | 0.2 |
| 14 | 8.94e-04 | 1.0 | 7.76e-05 | 1.0 | 2.50e-01 | -0.0 | 1.31e-01 | -1.3 | 1.25e-01 | 0.2 |

Figure 2 shows the computed densities $\rho$ at time $t = 1.8$ for $N = 800$ grid cells in $[-5, 5]$ ($h = 0.0125$) and time step size $\tau = 0.225h$. The shock front is well resolved in all cases. Regarding the fine structures in the density profile, parameter set $\Sigma_1$ is superior in this example, which matches well with the observations in [20]. The important point for the present paper is, however, that the differences between the results with and without artificial boundary at $x = 0$ are rather small for both $\Sigma_1$ and $\Sigma_2$.

### 4.3 Traffic network

We consider the road network depicted in Figure 3 consisting of three roads, a merging node $M$ and a dispersing node $D$.

On each road (with length $L = 0.2$) we consider the scalar conservation law

$$\rho_t + f(\rho)_x = 0$$

with flux function

$$f(\rho) = \begin{cases} \rho v & \text{if } 0 \le \rho \le \frac{\rho_{\max}}{2}, \\ (\rho_{\max} - \rho)\, v & \text{if } \frac{\rho_{\max}}{2} \le \rho \le \rho_{\max}, \end{cases}$$

where we use $v = 1$ und $\rho_{\max} = 1$.

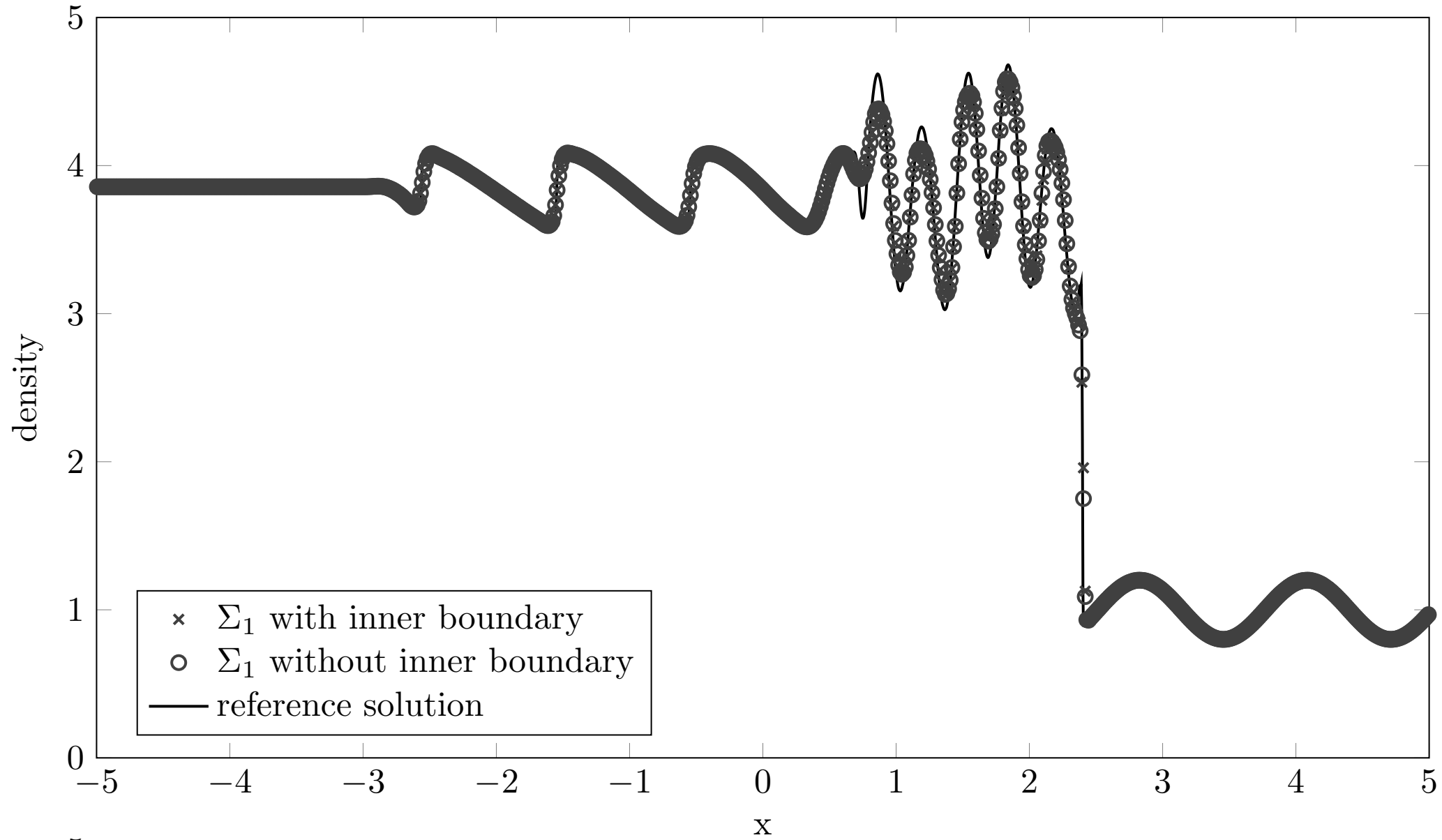

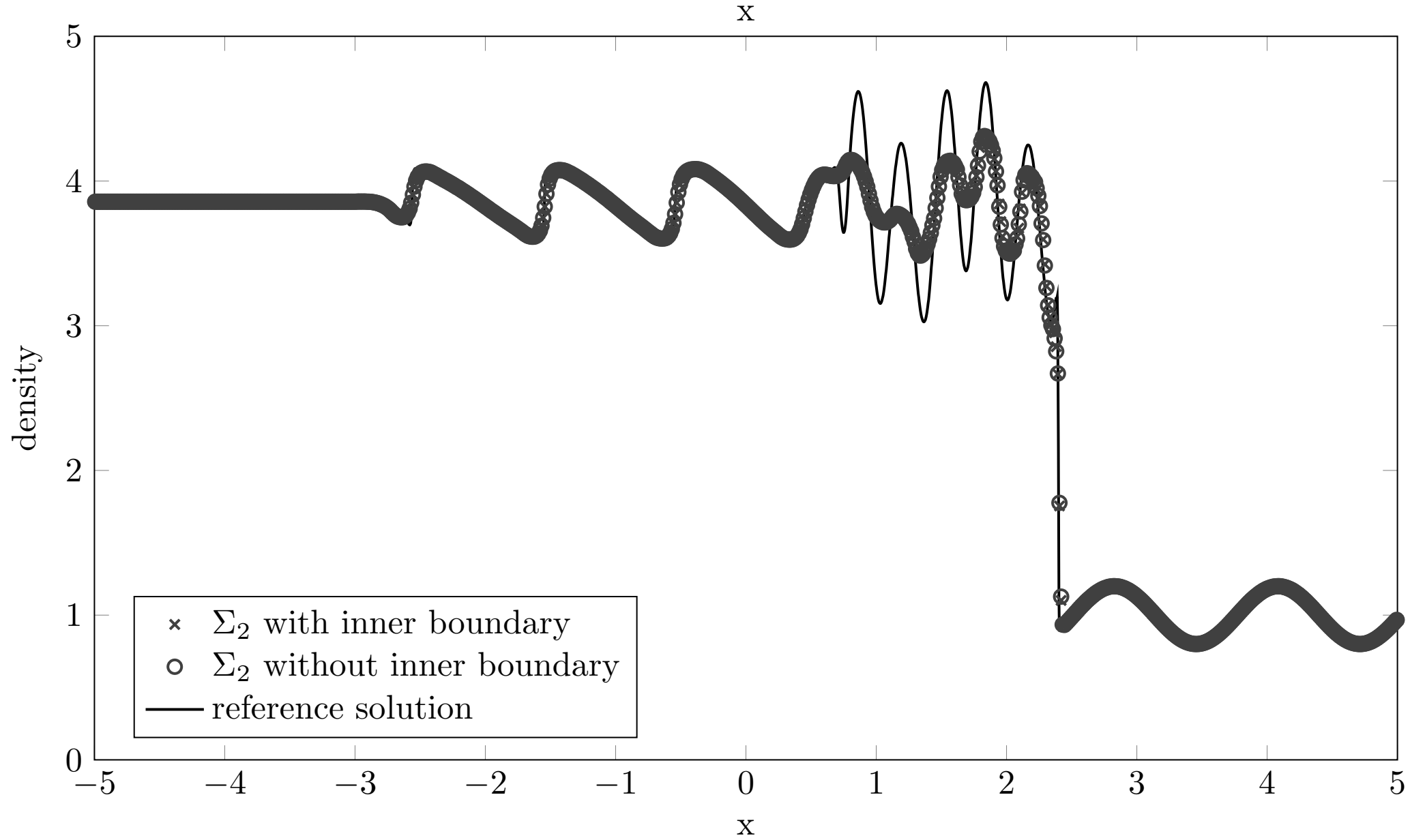

Figure 2: Results for the shock-acoustic interaction example with and without artificial inner boundary and with parameter settings $\Sigma_1$ and $\Sigma_2$. $N = 800$ cells are employed in this test.

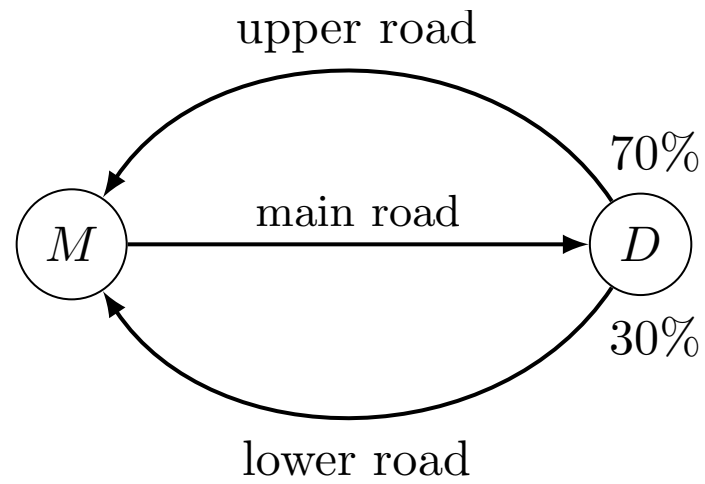

Figure 3: Road network with a merging and a dispersing junction.

At the diverging junction $D$ 70% of the flux go onto the upper road and 30% onto the lower road. The node $M$ is a merging junction with priority parameters (summing up to 1) for both incoming roads. For details on the applied junction models, we refer to [15].

First, we consider an accuracy test where we use priority parameter $P = 0.5$ for both incoming roads at the node $M$, and the following initial conditions:

$$\begin{aligned}
\text{main road:} \qquad & \rho_0(x) = 0.20 + 0.150 \sin\left(\pi \frac{x}{L} - \frac{\sin(\pi \frac{x}{L})}{\pi}\right), \\
\text{upper road:} \qquad & \rho_0(x) = 0.14 + 0.105 \sin\left(\pi \frac{x+L}{L} - \frac{\sin(\pi \frac{x+L}{L})}{\pi}\right), \\
\text{lower road:} \qquad & \rho_0(x) = 0.06 + 0.045 \sin\left(\pi \frac{x+L}{L} - \frac{\sin(\pi \frac{x+L}{L})}{\pi}\right).
\end{aligned}$$

Since we are in the free flow regime here ($\rho(x,t) \leq \frac{\rho_{\max}}{2}$ on all roads all the time), the analytical solution is known and equals the initial conditions at time $t = \frac{2L}{v} = 0.4$. For the numerical simulation we apply the parameter sets $\Sigma_1$ and $\Sigma_2$ and consider $N = 10 \cdot 2^n$ grid cells on each road ($n \in \{0, \ldots, 7\}$), i.e., $h = 0.02 \cdot 2^{-n}$, and $\tau = \frac{h}{2}$ as time step size.

Table 4 shows the $L^\infty$ error over the whole network at the final time $t = 0.4$ and the corresponding convergence rates. Obviously, third order convergence is achieved for the entire scheme on the whole network with both parameter sets.

Table 4: Simulation results for the traffic network depicted in Figure 3 with smooth initial conditions.

| | $n$ | 0 | 1 | 2 | 3 | 4 | 5 | 6 | 7 |
|---|---|---|---|---|---|---|---|---|---|
| $\Sigma_1$ | Error | 6.44e-03 | 9.37e-04 | 1.39e-04 | 1.90e-05 | 2.47e-06 | 3.14e-07 | 3.96e-08 | 4.97e-09 |
| | c.r. | - | 2.8 | 2.8 | 2.9 | 2.9 | 3.0 | 3.0 | 3.0 |
| $\Sigma_2$ | Error | 1.71e-02 | 4.42e-03 | 8.18e-04 | 1.17e-04 | 1.51e-05 | 1.88e-06 | 2.34e-07 | 2.92e-08 |
| | c.r. | - | 2.0 | 2.4 | 2.8 | 3.0 | 3.0 | 3.0 | 3.0 |

In the second test case, we consider a traffic jam on the upper road, which propagates backwards through the node $D$ onto the main road. For this, we consider a low priority parameter of $P = 0.2$ for the upper road at the merging node $M$, and the following initial conditions:

$$\begin{aligned}
\text{main road:} \qquad & \rho_0(x) = 0.50, \\
\text{upper road:} \qquad & \rho_0(x) = \begin{cases} 0.35 & \text{for } 0 \leq x \leq 0.05, \\ 0.90 & \text{for } 0.05 < x \leq 0.2, \end{cases} \\
\text{lower road:} \qquad & \rho_0(x) = \begin{cases} 0.15 & \text{for } 0 \leq x \leq 0.05, \\ 0.85 & \text{for } 0.05 < x \leq 0.2. \end{cases}
\end{aligned}$$

Note that the initial traffic jam on the lower road begins to resolve in the course of the simulation due to the higher priority parameter ($1.0 - P = 0.8$) of this road.

Simulation results for this traffic jam scenario with $N = 80$ cells per road ($h = 0.0025$), $\tau = \frac{h}{2}$ as time step size and parameter settings $\Sigma_1$ and $\Sigma_2$ are shown in Figure 4. The left plot shows the back-travelling discontinuity on the upper road at time $t = 0.05$, where almost no oscillations can be observed for both parameter sets. The plot on the right-hand side shows the jam front on the main road at time $t = 0.15$. Obviously, the results for parameter set $\Sigma_2$ are still non-oscillatory whereas the results for $\Sigma_1$ show a slight overshoot.

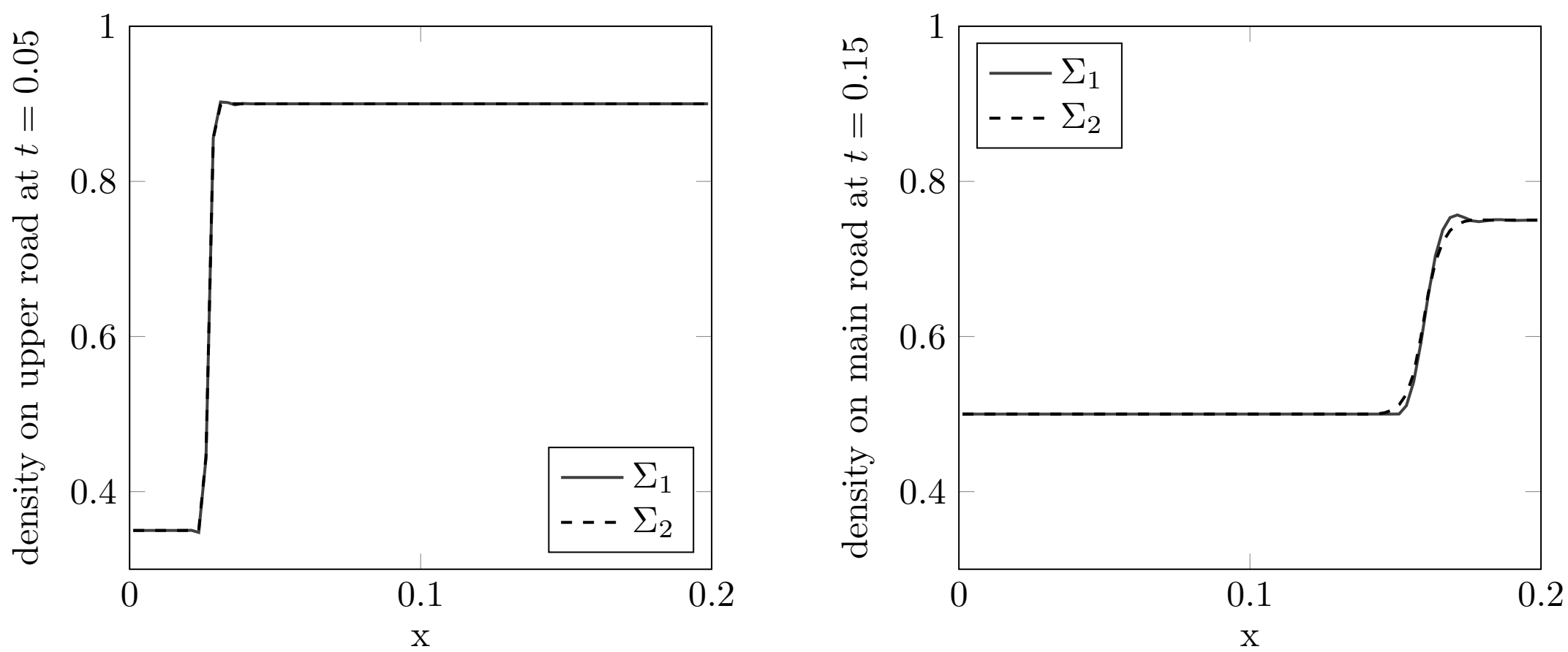

Figure 4: Results for the network example with traffic jam with parameter settings $\Sigma_1$ and $\Sigma_2$. $N = 80$ cells per road ($h = 0.0025$) are employed in this test.

## 4.4 Time-dependent boundary conditions

Next, we take a look at time-dependent boundary conditions. Therefore, we consider the traffic model from Section 4.3 for a single road with length $L = 0.4$, $v = 1$, $\rho_{\max} = 1$ and initial conditions

$$\rho_0(x) = 0.20 + 0.150 \sin\left(2\pi \frac{x}{L} - \frac{\sin(2\pi \frac{x}{L})}{\pi}\right).$$

On the left boundary we prescribe

$$\rho(0,t) = 0.20 + 0.150 \sin\left(-2\pi \frac{t}{L} - \frac{\sin(-2\pi \frac{t}{L})}{\pi}\right),$$

which reproduces the initial conditions at time $t = \frac{L}{v} = 0.4$. Numerically, the given boundary conditions directly lead to a time-dependent inflow condition (for all stages in the applied Runge-Kutta method). Note that such a straightforward treatment may give rise to order reduction effects as analysed in [7].

Table 5 shows the $L^\infty$ error at the final time $t = 0.4$ and the corresponding convergence rates (for $h = 0.02 \cdot 2^{-n}$ and $\tau = \frac{h}{2}$). For both parameter sets, full third order can be observed.

## 4.5 Open channel flow

We consider the network of channels depicted in Figure 5 with lengths $|\overline{AB}| = 0.1$, $|\overline{BC}| = 0.05$, $|\overline{BD}| = 0.05$, $|\overline{BE}| = 0.15$, $|\overline{CD}| = 0.05$, $|\overline{DE}| = 0.05$ and $|\overline{EA}| = 0.05$.

Table 5: Simulation results for a single road with time-dependent boundary conditions.

| $n$ | | 0 | 1 | 2 | 3 | 4 | 5 | 6 | 7 |
|---|---|---|---|---|---|---|---|---|---|
| $\Sigma_1$ | Error | 4.54e-03 | 6.88e-04 | 9.13e-05 | 1.17e-05 | 1.48e-06 | 1.86e-07 | 2.34e-08 | 2.92e-09 |
| | c.r. | - | 2.7 | 2.9 | 3.0 | 3.0 | 3.0 | 3.0 | 3.0 |
| $\Sigma_2$ | Error | 1.53e-02 | 3.75e-03 | 6.77e-04 | 1.02e-04 | 1.37e-05 | 1.75e-06 | 2.21e-07 | 2.78e-08 |
| | c.r. | - | 2.0 | 2.5 | 2.7 | 2.9 | 3.0 | 3.0 | 3.0 |

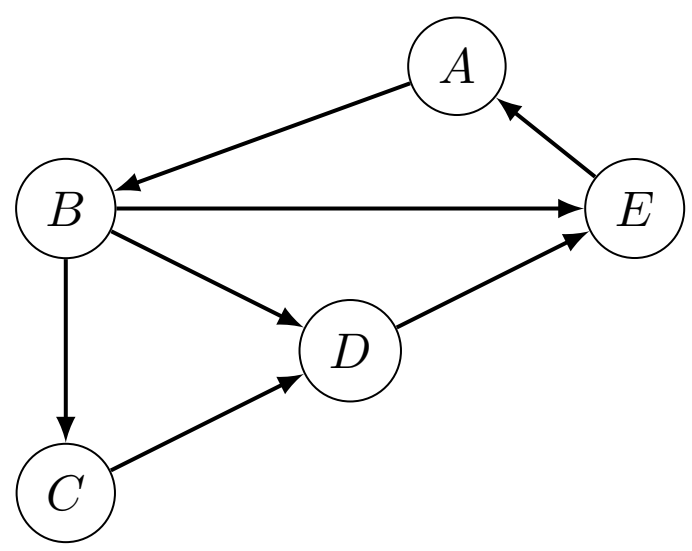

Figure 5: Network of open channels.

In each channel we consider the shallow water equations,

$$\frac{\partial}{\partial t}\begin{pmatrix} h \\ q \end{pmatrix} + \frac{\partial}{\partial x}\begin{pmatrix} q \\ \frac{q^2}{h} + \frac{1}{2}gh^2 \end{pmatrix} = 0 \tag{14}$$

with $g = 1.0$ here. At the intersections we prescribe conservation of mass (sum of ingoing fluxes equals sum of outgoing fluxes), continuity of the height $h$, and use the information of the ingoing characteristics of each channel. Since we consider an example in the subcritical regime, this procedure yields exactly the necessary number of equations at each intersection. Initially, we set the flux $q$ equal to zero in all channels, $h(x,0) = 0.3 + 0.03 \cdot (\sin(\pi \frac{x}{0.1}))^4$ in the channel from $A$ to $B$, and $h(x,0) = 0.3$ within all other channels.

For the numerical simulation we apply the parameter sets $\Sigma_1$ and $\Sigma_2$ and consider a spatial grid size $h = 0.01 \cdot 2^{-n}$ for all channels ($n \in \{0,\dots,7\}$) and $\tau = \frac{h}{2}$ as time step size until the final time $t = 0.2$. A reference solution is computed with the same parameters and $n = 9$. Table 6 shows the $L^\infty$ error over all states on the whole network at the final time and the corresponding convergence rates. Again, third order convergence is achieved for the entire scheme on the whole network with both parameter sets.

Table 6: Simulation results for the channel network depicted in Figure 5.

| $n$ | | 0 | 1 | 2 | 3 | 4 | 5 | 6 | 7 |
|---|---|---|---|---|---|---|---|---|---|
| $\Sigma_1$ | Error | 8.68e-04 | 2.78e-04 | 6.03e-05 | 1.15e-05 | 1.67e-06 | 2.20e-07 | 2.82e-08 | 3.51e-09 |
| | c.r. | - | 1.6 | 2.2 | 2.4 | 2.8 | 2.9 | 3.0 | 3.0 |
| $\Sigma_2$ | Error | 8.60e-04 | 8.52e-04 | 6.59e-05 | 1.11e-05 | 1.53e-06 | 2.01e-07 | 2.56e-08 | 3.18e-09 |
| | c.r. | - | 0.1 | 3.7 | 2.6 | 2.8 | 2.9 | 3.0 | 3.0 |

## 4.6 Dam break problem

This test is designed to test the performance of the reconstruction procedure when a shock reflects from a static wall. We employ again the shallow water equations (14) in the domain $x \in [0,1]$, but employ wall boundary conditions and set up dam-break initial conditions: the initial velocity is zero everywhere and the initial water height is

$$(a): \ h_0(x) = \begin{cases} 1, & x \in [0.4, 0.6], \\ 0.5, & x \notin [0.4, 0.6], \end{cases} \qquad (b): \ h_0(x) = \begin{cases} 1, & x \in [0.4, 0.6], \\ 0.05, & x \notin [0.4, 0.6]. \end{cases}$$

In both cases two shocks move towards the sides of the domain and two rarefactions towards the center. In case $(a)$ the flow is always subcritical (Froude number less than 0.36), while the second one has supercritical regions (with Froude numbers up to 1.72). The final time is chosen as $t = 0.6$, so that the shocks have just reflected from the walls at the ends of the domain.

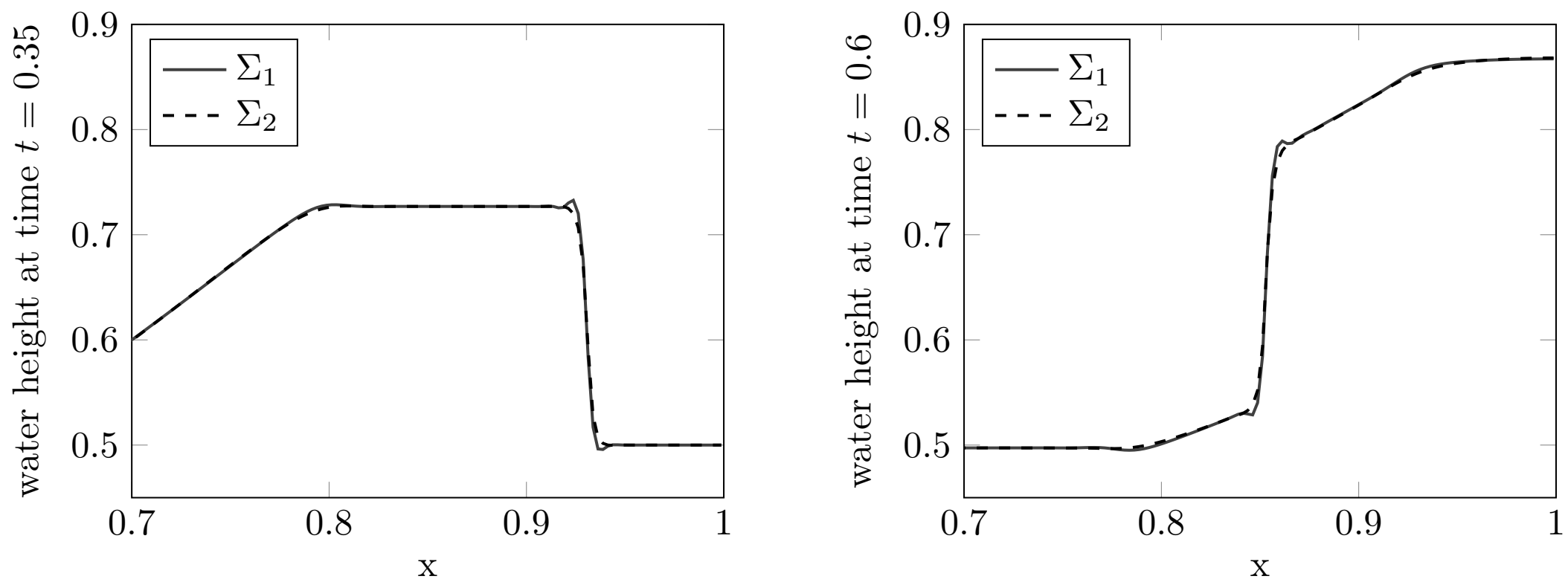

Figure 6: Dam break scenario with wall boundary conditions. Results for initial data $(a)$, $N = 400$ grid cells ($h = 0.0025$) at times $t = 0.35$ (left) and $t = 0.6$ (right).

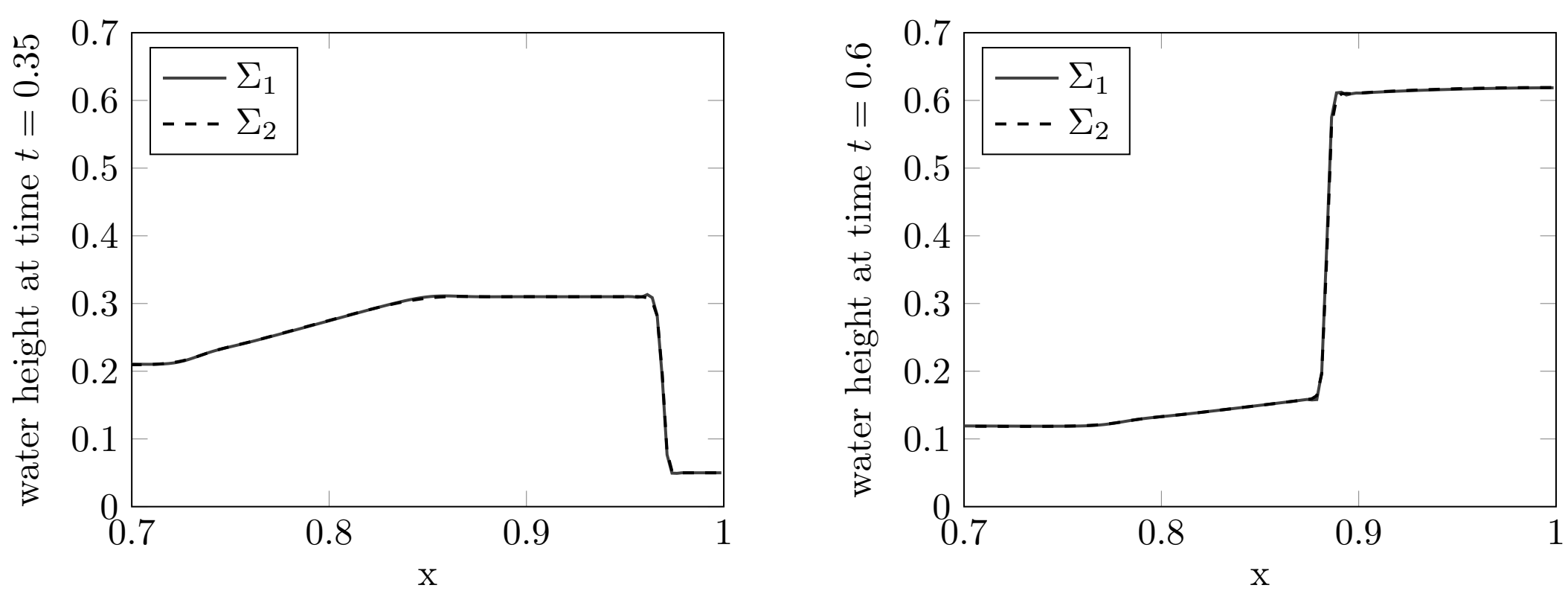

Figure 7: Dam break scenario with wall boundary conditions. Results for initial data $(b)$, $N = 400$ grid cells ($h = 0.0025$) at times $t = 0.35$ (left) and $t = 0.6$ (right).

Note that, when the shocks are about to hit the wall, there is no smooth upwind stencil to be chosen and therefore a reconstruction procedure that wasn't able to drop down to a constant reconstruction would introduce oscillations due to the employment of downwind and/or oscillatory stencils. We present the results obtained with the proposed reconstruction for parameter sets $\Sigma_1$ and $\Sigma_2$ in Figure 6 for initial data $(a)$ and in Figure 7 for initial data $(b)$, both at times $t = 0.35$ (before the shock reflection) and $t = 0.6$ (after the shock reflection), where $N = 400$ grid cells have been used ($h = 0.0025$) and time step size $\tau = 0.45h/\lambda$ with $\lambda = \max |q/h| + \sqrt{gh}$ in each time step. Similar to the results in Section 4.3 for the traffic jam situation, the results for $\Sigma_2$ are completely non-oscillatory whereas the results for $\Sigma_1$ show slight oscillations. Note that the observed over-/undershoots for $\Sigma_1$ are already present in the numerical solution before the shock reflection at the wall and are not introduced by the boundary treatment.

# 5 Conclusion

We considered the convergence properties of a third order CWENO boundary treatment. Based on the ideas of the WENO extrapolation of [30], we proved sufficient conditions for the underlying scheme parameters to achieve the optimal convergence rates in different cases. All theoretical findings could be confirmed by the presented numerical results. In particular, full third order could be observed for a smooth test example on a small traffic network and a network of channels governed by the shallow water equations. Further, a non-oscillatory behaviour of the proposed boundary treatment could be observed for discontinuities passing a boundary within the same traffic network, for the well-known shock-acoustic interaction example and for shock reflections from a wall for the shallow water equations.

As mentioned in Remark 2.2, we point out that the well-known linear scaling limiter from [31, 32] can be directly integrated in the presented boundary treatment as well as in the applied third order central WENO scheme in the interior of the computational domain. This way one may directly achieve a maximum-principle-satisfying and positivity-preserving third order scheme on the whole network.

In future work, we intend to generalize and analyse the proposed boundary treatment in the case of non-uniform grids and for higher dimensions. Further, generalizations of the results to higher than third order schemes are within our scope as well as other properties relevant for specific cases like well-balancing.

# References

[1] Aràndiga, F., Baeza, A., Belda, A. M., and Mulet, P. Analysis of WENO Schemes for Full and Global Accuracy. *SIAM J. Numer. Anal. 49*, 2 (2011), 893–915.

[2] Banda, M. K., Häck, A.-S., and Herty, M. Numerical Discretization of Coupling Conditions by High-Order Schemes. *J. Sci. Comput. 69*, 1 (2016), 122–145.

[3] Borsche, R., and Kall, J. ADER schemes and high order coupling on networks of hyperbolic conservation laws. *J. Comput. Phys. 273* (2014), 658 – 670.

[4] Borsche, R., and Kall, J. High order numerical methods for networks of hyperbolic conservation laws coupled with ODEs and lumped parameter models. *J. Comput. Phys. 327* (2016), 678 – 699.

[5] Borsche, R., Klar, A., Kühn, S., and Meurer, A. Coupling Traffic Flow Networks to Pedestrian Motion. *Math. Models Methods Appl. Sci. 24*, 02 (2014), 359–380.

[6] Brouwer, J., Gasser, I., and Herty, M. Gas pipeline models revisited: Model hierarchies, nonisothermal models, and simulations of networks. *Multiscale Model. Simul. 9*, 2 (2011), 601–623.

[7] Carpenter, M. H., Gottlieb, D., Abarbanel, S., and Don, W.-S. The theoretical accuracy of runge–kutta time discretizations for the initial boundary value problem: A study of the boundary error. *SIAM J. Sci. Comput. 16*, 6 (1995), 1241–1252.

[8] Castro, M., Costa, B., and Don, W. S. High order weighted essentially non-oscillatory WENO-Z schemes for hyperbolic conservation laws. *J. Comput. Phys. 230*, 5 (2011), 1766–1792.

[9] Coclite, G. M., Garavello, M., and Piccoli, B. Traffic Flow on a Road Network. *SIAM J. Math. Anal. 36*, 6 (2005), 1862–1886.

[10] Contarino, C., Toro, E. F., Montecinos, G. I., Borsche, R., and Kall, J. Junction-Generalized Riemann Problem for stiff hyperbolic balance laws in networks: An implicit solver and ADER schemes. *J. Comput. Phys. 315* (2016), 409 – 433.

[11] Cravero, I., and Semplice, M. On the Accuracy of WENO and CWENO Reconstructions of Third Order on Nonuniform Meshes. *J. Sci. Comput. 67*, 3 (2016), 1219–1246.

[12] Domschke, P., Kolb, O., and Lang, J. Adjoint-based error control for the simulation and optimization of gas and water supply networks. *Appl. Math. Comput. 259* (2015), 1003–1018.

[13] Don, W.-S., and Borges, R. Accuracy of the weighted essentially non-oscillatory conservative finite difference schemes. *J. Comput. Phys. 250* (2013), 347–372.

[14] Feng, H., Hu, F., and Wang, R. A new mapped weighted essentially non-oscillatory scheme. *J. Sci. Comput. 51*, 2 (2012), 449–473.

[15] Goatin, P., Göttlich, S., and Kolb, O. Speed limit and ramp meter control for traffic flow networks. *Eng. Optim. 48*, 7 (2016), 1121–1144.

[16] Gottlieb, S., and Shu, C.-W. Total variation diminishing Runge-Kutta schemes. *Math. Comp. 67* (1998), 73–85.

[17] Göttlich, S., Kolb, O., and Kühn, S. Optimization for a special class of traffic flow models: Combinatorial and continuous approaches. *Netw. Heterog. Media 9*, 2 (2014), 315–334.

[18] Ha, Y., Ho Kim, C., Ju Lee, Y., and Yoon, J. An improved weighted essentially non-oscillatory scheme with a new smoothness indicator. *J. Comput. Phys. 232*, 1 (2013), 68–86.

[19] Jiang, G.-S., and Shu, C.-W. Efficient Implementation of Weighted ENO Schemes. *J. Comput. Phys. 126*, 1 (1996), 202–228.

[20] Kolb, O. On the Full and Global Accuracy of a Compact Third Order WENO Scheme. *SIAM J. Numer. Anal. 52*, 5 (2014), 2335–2355.

[21] Kolb, O. On the Full and Global Accuracy of a Compact Third Order WENO Scheme: Part II. In *Numerical Mathematics and Advanced Applications ENUMATH 2015*, B. Karasözen, M. Manguoğlu, M. Tezer-Sezgin, S. Göktepe, and Ö. Uğur, Eds. Springer International Publishing, 2016, pp. 53–62.

[22] Kolb, O. A Third Order Hierarchical Basis WENO Interpolation for Sparse Grids with Application to Conservation Laws with Uncertain Data. *J. Sci. Comput.* (2017).

[23] Kolb, O., and Lang, J. *Mathematical optimization of water networks*, vol. 162 of *Internat. Ser. Numer. Math.* Birkhäuser/Springer Basel AG, 2012, ch. Simulation and continuous optimization, pp. 17–33.

[24] Levy, D., Puppo, G., and Russo, G. Compact Central WENO Schemes for Multidimensional Conservation Laws. *SIAM J. Sci. Comput. 22*, 2 (2000), 656–672.

[25] Liu, X.-D., Osher, S., and Chan, T. Weighted Essentially Non-oscillatory Schemes. *J. Comput. Phys. 115*, 1 (1994), 200–212.

[26] Müller, L. O., and Blanco, P. J. A high order approximation of hyperbolic conservation laws in networks: Application to one-dimensional blood flow. *J. Comput. Phys. 300* (2015), 423 – 437.

[27] Müller, L. O., and Toro, E. F. A global multiscale mathematical model for the human circulation with emphasis on the venous system. *Int. J. Numer. Methods Biomed. Eng. 30*, 7 (2014), 681–725.

[28] Semplice, M., Coco, A., and Russo, G. Adaptive Mesh Refinement for Hyperbolic Systems Based on Third-Order Compact WENO Reconstruction. *J. Sci. Comput. 66*, 2 (2016), 692–724.

[29] Shu, C.-W., and Osher, S. Efficient implementation of essentially non-oscillatory shock-capturing schemes, II. *J. Comput. Phys. 83*, 1 (1989), 32–78.

[30] Tan, S., and Shu, C.-W. Inverse Lax-Wendroff procedure for numerical boundary conditions of conservation laws. *J. Comput. Phys. 229*, 21 (2010), 8144–8166.

[31] Zhang, X., and Shu, C.-W. On maximum-principle-satisfying high order schemes for scalar conservation laws. *J. Comput. Phys. 229*, 9 (2010), 3091–3120.

[32] Zhang, X., and Shu, C.-W. Maximum-principle-satisfying and positivity-preserving high-order schemes for conservation laws: survey and new developments. *Proc. R. Soc. Lond. Ser. A Math. Phys. Eng. Sci. 467*, 2134 (2011), 2752–2776.